\documentclass{amsart}

\usepackage{amssymb}
\usepackage[all]{xy}
\usepackage{hyperref}

\usepackage{enumitem}   

\setlist[enumerate]{itemsep=.2em,topsep=.2em,leftmargin=1.25em,itemindent=2.0em}

\newtheorem{thm}{Theorem}

\newtheorem{lem}[thm]{Lemma}
\newtheorem{cor}[thm]{Corollary}

\newtheorem{prop}[thm]{Proposition}

\theoremstyle{definition}
\newtheorem{defn}[thm]{Definition}

\newtheorem{say}[thm]{}
\newtheorem{exmp}[thm]{Example}

\newtheorem{rem}[thm]{Remark}          
\newtheorem{aside}[thm]{Aside}          
\newtheorem*{ack}{Acknowledgments}      

\newtheorem{defn-thm}[thm]{Definition--Theorem}  
\newtheorem{defn-lem}[thm]{Definition--Lemma}  

\newtheorem{assertion}[thm]{Assertion}   
   
\newtheorem{vardefs}[thm]{Variant Definitions}   
\newtheorem{main-exmp}[thm]{Main Example}
\newtheorem{baby-exmp}[thm]{Baby Example}
\newtheorem*{obs-nn}{Observation}
\newtheorem*{assertion-nn}{Assertion}   

\theoremstyle{remark}

\let \cedilla =\c
\renewcommand{\c}[0]{{\mathbb C}}  

\renewcommand{\o}[0]{{\mathcal O}} 
\newcommand{\z}[0]{{\mathbb Z}}
\newcommand{\n}[0]{{\mathbb N}}
\renewcommand{\r}[0]{{\mathbb R}} 

\renewcommand{\a}[0]{{\mathbb A}}

\newcommand{\p}[0]{{\mathbb P}}
\newcommand{\f}[0]{{\mathbb F}}

\newcommand{\map}[0]{\dasharrow}
\newcommand{\qtq}[1]{\quad\mbox{#1}\quad}
\newcommand{\spec}[0]{\operatorname{Spec}}
\newcommand{\pic}[0]{\operatorname{Pic}}
\newcommand{\pics}[0]{\operatorname{\mathbf{Pic}}}

\newcommand{\pictw}[0]{\operatorname{Pic}^{\rm tw}}

\newcommand{\gal}[0]{\operatorname{Gal}}

\newcommand{\rank}[0]{\operatorname{rank}}
\newcommand{\mult}[0]{\operatorname{mult}}

\newcommand{\proj}[0]{\operatorname{Proj}}

\newcommand{\coker}[0]{\operatorname{coker}}    
\newcommand{\ext}[0]{\operatorname{Ext}}    
\newcommand{\Hom}[0]{\operatorname{Hom}}    
\newcommand{\trace}[0]{\operatorname{trace}}

\newcommand{\aut}[0]{\operatorname{Aut}}

\newcommand{\chr}[0]{\operatorname{char}}

\newcommand{\br}[0]{\operatorname{Br}}

\newcommand{\onto}[0]{\twoheadrightarrow}

\newcommand{\tsum}[0]{\textstyle{\sum}}

\newcommand{\bir}[0]{\operatorname{Bir}} 
 
\newcommand{\ind}[0]{\operatorname{index}}

\def\into{\DOTSB\lhook\joinrel\to}

\def\loccoh#1.#2.#3.#4.{H^{#1}_{#2}(#3,#4)}

\DeclareMathAlphabet{\mathchanc}{OT1}{pzc}%
                                {m}{it}

\newcommand{\gm}[0]{{\mathbb G}_m}

\newcommand{\simb}[0]{\stackrel{bir}{\sim}}

\newcommand{\OG}{\mathrm{OG}}

\newcommand{\sym}[0]{\operatorname{Sym}}

\usepackage[all]{xy}\xyoption{dvips}

\newcommand{\tprod}[0]{\textstyle{\prod}} 

\newcommand{\norm}[0]{\operatorname{norm}}

\newcommand{\End}[0]{\operatorname{End}}
\newcommand{\sEnd}[0]{{\mathchanc{End}}}

\begin{document}
\bibliographystyle{amsalpha}

\today

 \title[Severi--Brauer varieties]{Severi--Brauer varieties;\\ a geometric treatment}
 \author{J\'anos Koll\'ar}

 \begin{abstract} We present a geometric treatment of Severi--Brauer varieties, without using any results from the theory of central simple algebras or   Galois cohomology.
   
       \end{abstract}

 \maketitle

\tableofcontents

The story of Severi--Brauer varieties ties together three---seemingly unrelated---strands.

$\bullet$ {\it Central simple algebras}  are finite dimensional, simple algebras, whose center is a given field $k$. The first example is the quaternions, discovered by  Hamilton   in 1843.    Wedderburn's 
contributions started  in  1905, and the subject grew into its own  around 
1925 with the works of Albert, Brauer, Hasse and  Noether.
See  \cite[Chap.13]{vdW-book}, \cite[Chap.4]{MR571884},
\cite[Chap.XVII]{lang-alg}, \cite{MR2266528}
or \cite[Chap.073W]{stacks-project}
for introductions.

$\bullet$ {\it Severi--Brauer varieties} are varieties defined over a field $k$, that become isomorphic to projective space over the algebraic closure $\bar k$. Plane conics without $k$-points give the `baby examples,' though they were usually viewed as part of the theory of quadratic forms.
 Higher dimensional Severi--Brauer varieties
 were first defined and studied by Severi \cite{severi-SB}, and their theory came into full maturity starting with  Ch\^atelet  \cite{MR0014720}, who used both  Galois cohomology and central simple algebras. 
For a variety $X$, 
maps to non-trivial Severi--Brauer varieties give the best understood  obstruction to having a $k$-point. Thus  Severi--Brauer varieties appear  especially frequently in   arithmetic questions, see  \cite{MR4304038}.

$\bullet$  $H^2(k, \gm)$, the second Galois cohomology group of the multiplicative group $\gm$.   This is the most powerful viewpoint
\cite{MR2266528}.
\medskip

The aim of these notes is to present a geometric treatment of Severi--Brauer varieties, without using any results from the theory of central simple algebras or   Galois cohomology. This reverses the 
  treatments prevalent since \cite{MR0014720}, which  start either with
 central simple algebras or with Galois cohomology, and proceed to derive geometric properties.

A  series of {\it Asides} outline  either
  connections to the other approaches to Severi--Brauer varieties, or discuss connections with other topics.
These are not needed for the proofs, but provide  motivation and hints for those who intend to go deeper into some of these topics. 

Introductions to the algebraic and cohomological methods are in 
\cite{MR2266528} and \cite[Chap.IV]{milne}.
The arithmetic aspects are treated in \cite{MR4304038}.
Connections with algebraic geometry are discussed in \cite{gro-bra-I-III, MR3587845}, while 
\cite{MR611862, MR657419} are more algebraic.

\medskip
The main idea of our treatment is the following.

\begin{assertion}\label{ass.1.1}  The geometry of  a  Severi--Brauer variety $P$  is  best studied via the vector bundle 
 $F(P)$ which is obtained as  the unique (up-to isomorphism)  non-split extension
$$
0\to \o_P\to F(P)\to T_P\to 0,
\eqno{(\ref{ass.1.1}.1)}
$$
where $T_P$ denotes the tangent bundle of $P$.
\end{assertion}

The vector bundle $F(P)$ appears in Quillen's work \cite[\S 8.4]{MR0338129}
(with a different construction),  a similar bundle is used in
an  unpublished preprint  of  Szab\'o
\cite{sz-s-b} (which  also contains the first steps of a  program to treat 
Severi--Brauer varieties geometrically), and is a special case of the
absolutely split bundles studied by
 Biswas--Nagaraj \cite{MR2391341} and 
Novakovi\'c  \cite{MR4796081}.

\medskip
{\bf Description of the Sections}
\medskip

We start in Section~\ref{sec.1} by studying the key abstract property of the vector bundle $F(P)$: it is a direct sum of copies of $\o(1)$  over $\bar k$.

The first  examples of  Severi--Brauer varieties are exhibited  in Section~\ref{sec.1.b}. Basic properties of twisted linear subvareties are established in
Section~\ref{sec.2} and splitting fields are constructed in
Section~\ref{sec.split.f}.

We show in   Section~\ref{sec.2.c} that the  examples of Section~\ref{sec.1.b}
give all  non-minimal Severi--Brauer varieties, and   that
$\aut(P)$ acts transitively on flags of  twisted linear subvareties.

Twisted linear systems are studied in Section~\ref{sec.6}, leading to a geometric definition of the  Brauer group in Section~\ref{sec.7}.
The index and the period of  Severi--Brauer varieties is studied in Section~\ref{sec.8}.

In  Section~\ref{sec.norm.form}  we construct
 cyclic  Severi--Brauer varieties  from hypersurfaces
of degree $n+1$ in $\p^{n+1}$ and norm forms.

In Section~\ref{sec.11} we construct
the sequence (\ref{ass.1.1}.1) for  Severi--Brauer schemes.

\begin{ack} I thank  J.-L.~Colliot-Th\'el\`ene,  A.J.~de~Jong, P.~Gille, D.~Litt, S.~ Novakovi\'c, M.~Shin, C.~Skinner, J.~Starr  and 
T.~Szamuely  for many
comments, corrections  and references.
 Partial  financial support    was provided  by  the NSF (grant number
DMS-1901855)  and by the Simons Foundation   (grant number SFI-MPS-MOV-00006719-02).
\end{ack}

\section{Geometrically split vector bundles}\label{sec.1}

In this section we study abstract properties of the vector bundle $F(P)$ defined in (\ref{ass.1.1}.1). Its key  feature is that,  while
$F(P)$ is usually indecomposable over the ground field $k$, it is a direct sum of several copies of a line bundle over $\bar k$. 
 Our main interest is in smooth projective varieties, but
the basic results in this section hold in the following setting.

\begin{defn} \label{T.L.defn}
Let $k$ be a field,  and  $\bar k$ a fixed  algebraic closure of $k$.
Let $X$ be a  proper  $k$-scheme that is geometrically reduced and  geometrically connected. (That is, $X_{\bar k}$ is reduced and connected.)

The basic results all work if $X$ is not proper, as long as $H^0(X,\o_X)\cong k$. 

Given 
a line bundle $ {\mathcal L}$  on $X_{\bar k}$,
 let $T( {\mathcal L})$ denote  the category of vector bundles $F$ on $X$ that 
$F_{\bar k}\cong \oplus_i {\mathcal L}$; a sum of copies of $ {\mathcal L}$.

${\mathcal L}$ is called a {\it twisted line bundle} on $X$  if 
$T({\mathcal L})$ contains a 
nonzero vector bundle.

Let ${\mathcal L}_i$ be twisted line bundles on $X$  and 
$0\neq F_i\in T({\mathcal L}_i)$.
Then $F_1\otimes F_2\in T({\mathcal L}_1\otimes {\mathcal L}_2)$,
thus all twisted line bundles on $X$ form a group
$\pictw(X)$.

({\it Aside.} We note in (\ref{rem.on.tw.l}.4) that $
\pictw(X)=\pics(X)(k)$, the $k$-points of the Picard scheme of $X$.)

A line bundle  $L$  on $X$  is naturally identified with the
twisted line bundle ${\mathcal L}:=L_{\bar k}$. 
We see in (\ref{split.lem})  that in this case 
$T({\mathcal L})$ consists of sums of copies of $ L$.
Conversely, if $T({\mathcal L})$ contains a rank 1 bundle $L$, then
${\mathcal L}\cong L_{\bar k}$. 
Thus there is a natural injection
$$
\pic(X)\into \pictw(X).
\eqno{(\ref{T.L.defn}.1)}
$$
If $M$ is a line bundle on $X$, then $T( {\mathcal L}\otimes M_{\bar k})=T( {\mathcal L})\otimes M$.
Thus twisted line bundles really study the quotient 
$\pictw(X)/\pic(X)$.

Technically it is better to  work with a separable algebraic closure $ k^s$ instead of an algebraic closure  $\bar k$. (Note that $\bar k=k^s$ if $\chr k=0$, or if
$\chr k>0$ and $k$ is perfect.)
We prove in  (\ref{SB.triv.sc.cor})   
that the two versions are equivalent.

Our definition of  twisted line bundles is related to, but quite different from,
the notion of twisted coherent sheaves studied in \cite{MR2306170, MR2309155}.
\end{defn}

The simplest example of a twisted line bundle is $\o_C(1)$ of a plane conic. 

\begin{baby-exmp}  \label{babay.exmp}
Let $C$ be a $k$-variety such that $C_{\bar k}\cong \p^1_{\bar k}$.
Then $\o_C(-K_C)$ is very ample and gives an embedding
$C\into \p^2_k$ whose image is a conic; see \cite[Secs.1.2--3]{ksc} for details.

Assume that $C$ has no $k$-points. 
(For example, if $k\subset \r$ then we can take
$C=(x^2+y^2+z^2=0)$.)

Then $C$ has no degree 1 line bundles defined over $k$ (since the zero set of any section would be a $k$-point), but  we have
$\o_C$ and $\o_C(2)=\o_{\p^2}(1)|_C$. Furthermore, 
$H^1\bigl(C, \o_C(-2)\bigr)\cong k$, thus we have a non-split extension
$$
0\to \o_C\to F(C)\to \o_C(2)\to 0.
\eqno{(\ref{babay.exmp}.1)}
$$
We know that, over $\bar k$,   $F(C)$ is isomorphic to
$\o_{C_{\bar k}}(1)+\o_{C_{\bar k}}(1)$, but $F(C)$ itself is an indecomposable vector bundle over $k$.  Thus $F(C)\in T({\mathcal L})$
for ${\mathcal L}:=\o_{C_{\bar k}}(1)$.
We see next that it is better to think of (\ref{babay.exmp}.1) as
$$
0\to \o_C\to F(C)\to T_C\to 0.
\eqno{(\ref{babay.exmp}.2)}
$$
\end{baby-exmp}

\begin{defn}[Severi--Brauer variety]\label{SB.s.1.defn}
  A {\it Severi--Brauer variety} over a field $k$ is a
variety $P$ such that $P_{\bar k}\cong \p^n_{\bar k}$, where $\bar k$ is an algebraic closure of $k$, and $n=\dim P$.

We  see in  (\ref{SB.triv.sc.cor}) that this is equivalent
to  $P_{k^s}\cong \p^n_{k^s}$ where $k^s$ is a separable algebraic closure of $k$.

It is not easy  to write down explicit  Severi--Brauer varieties of dimension  $n\geq 2$ or even to show that they exist; see Sections~\ref{sec.1.b} and \ref{sec.norm.form} for some examples.

More generally, 
let $K/k$ be a field extension and $X$ a $K$-scheme.
A $k$-scheme $Y$ is called a $k$-form of $X$ if $X\cong Y\times_{\spec k}\spec K$; we  abbreviate this as $X\cong Y_K$.
 See \cite[Sec.V.20]{MR0103191}, 
\cite[Sec.V.20]{MR918564} or \cite[pp.456-463]{MR2675155}
 for  general discussions.

Thus  $n$-dimensional Severi--Brauer varieties over $k$ are  the $k$-forms of $\p^n_{\bar k}$. 
\end{defn}

\begin{main-exmp}[$\o_P(1)$ of a Severi--Brauer variety] \label{main.exmp}
Let $P$ be a Severi--Brauer variety over a field $k$.

As in the 1-dimensional case,  $\o_P(-K_P)$ is very ample and gives an embedding
$P\into \p^N_k$ where $N=\binom{2n+1}{n}-1$, but this is very hard to use  since the codimension of $P$ is  large.

The tangent bundle of $\p^n$ can be presented using the Euler sequence
$$
0\to \o_{\p^n}\stackrel{e}{\to} \oplus_{i=0}^n\o_{\p^n}(1)
\stackrel{\partial}{\to}  T_{\p^n}\to 0;
\eqno{(\ref{main.exmp}.1)}
$$
 see   \cite[3.2.4]{eh-3264} or 
 \cite[II.8.13]{hartsh}.

Note that the sequence (\ref{main.exmp}.1) is determined by
its nonzero extension class  
$$
\eta\in \ext^1(T_{\p^n}, \o_{\p^n})=H^1({\p^n}, \Omega_{\p^n}^1)\cong  k.
\eqno{(\ref{main.exmp}.2)}
$$
Thus, for any Severi--Brauer variety $P$,
$$
 \ext^1(T_{P}, \o_{P})=H^1(P, \Omega_P^1)\cong  k.
\eqno{(\ref{main.exmp}.3)}
$$
(It is not important for now, but in  (\ref{R1.lem})  we  write down a
canonical isomorphism $H^1(P, \Omega_P^1)\cong k$.)

We can thus define a vector bundle $F(P)$ as the unique (up-to scaling) non-split extension
$$
0\to \o_P\to F(P)\to T_P\to 0.
\eqno{(\ref{main.exmp}.4)}
$$
Note that $F(P)_{\bar k}\cong \oplus_{i=0}^n\o_{\p^n}(1)$, thus
$F(P)\in T({\mathcal L})$ where  ${\mathcal L}=\o_{\p^n_{\bar k}}(1)$.

The Picard group of
$\p^n$ is  $\z[\o_{\p^n}(1)]\cong \z$, but, as in the baby example,
there is usually no line bundle $ \o_P(1)$. 
In fact, we see during the first proof of (\ref{split.lem.cor.4}) that $ \o_P(1)$ exists as a line bundle on $P$ iff $P\cong \p^{\dim P}$.
(In most cases, the only line bundles on $P$ are given by (multiples of)
the canonical class $K_P$, equivalently, by  (powers of) the dualizing sheaf $\omega_P$. Both $\o_P(K_P)$ and $\omega_P$
correspond to   $\o_P(-\dim P-1)$.)

This gives us
the main examples  of varieties with twisted line bundles: $X=P$, ${\mathcal L}=\o_{P_{\bar k}}(1)$.  
We set 
$$
T(P):=T\bigl(\o_{P_{\bar k}}(1)\bigr) \qtq{and let} E(P):=E\bigl(\o_{P_{\bar k}}(1)\bigr)
\eqno{(\ref{main.exmp}.5)}
$$ 
denote a minimal rank direct summand of $F(P)$.
(We see in (\ref{split.lem}) that $E(P)$ is unique up to isomorphism.) 

It follows from 
(\ref{split.lem.cor.2})  and (\ref{tw.lb.defn}.2) that the examples  $\bigl(P, {\mathcal L}\bigr)$ are universal: every other 
twisted line bundle  is obtained from some $\bigl(P, {\mathcal L}\bigr)$ by pull-back and tensoring by a (non-twisted) line bundle.

{\it Note.} We chose to write everyting in terms of the 
twisted line bundle $\o_{P_{\bar k}}(1)$ since  the resulting vector bundle $F$ fits in the exact sequence (\ref{main.exmp}.4).  One could also use its dual
$\o_{P_{\bar k}}(-1)$. That would give $F^*$ which is related to the cotangent bundle $\Omega_P^1$. More generally, the other twisted line bundles  $\o_{P_{\bar k}}(m)$ also lead to  interesting vector bundles. If $(m, \dim P+1)=1$ then
$T\bigl(\o_{P_{\bar k}}(m)\bigr)$ is essentially equivalent to
$T(P)$.
\end{main-exmp}

\begin{vardefs} \label{vardefs.1}
Although we will only deal with $T({\mathcal L})$ as defined above, 
the definition  makes sense if $X$ is normal and ${\mathcal L}$ is a
rank 1 reflexive sheaf. All the properties in this note hold in that setting. This is sometimes convenient in studying  the pull-back of a line bundle by a rational map $f:X\map Y$.   If $Y$ is proper and $X$ is smooth then $f^*{\mathcal L}$ naturally  extends to a line bundle on $X_{\bar k}$, but if
$X$ is only normal then $f^*{\mathcal L}$  usually extends only to a
rank 1 reflexive sheaf.

The definition and the basic properties also work when
${\mathcal L}$ is replaced by a  {\it simple} coherent sheaf ${\mathcal G}$ on $X_{\bar k}$.  (Simple means that  $\End({\mathcal G})\cong \bar k$.) 

It is also of interest to investigate   bundles $F$ such that
$F_{\bar k}$ is a direct sum of line bundles of different degrees. 
The study of such bundles on  Severi--Brauer varieties  was started in \cite{MR2391341, MR2507589} and completed in \cite{MR4796081}, where they are called {\it absolutely split} bundles. These bundles play a key role in understanding the bounded derived category of coherent sheaves on  Severi--Brauer varieties; see \cite{MR2863422} and \cite{2017arXiv170103020N}.

{\it Aside \ref{vardefs.1}.1.} One can also consider those bundles $F$ such that
$F_{\bar k}$ is an {\em extension} of copies of ${\mathcal L}$, or even of line bundles numerically equivalent to ${\mathcal L}$.  This gives a better notion on Abelian varieties, though the resulting category is not semisimple; 
see \cite{MR607081, MR2529476, MR3194649} for various characterizations and uses of these vector bundles.
\end{vardefs}

\begin{aside}\label{rem.on.tw.l}  Let $X$ be any $k$-scheme, $F$ a coherent sheaf on $X$,  and $K/k$ a Galois extension.
Then $F_K$ is a   coherent sheaf on $X_K$ and 
$$
F_K^{\sigma}\cong F_K\qtq{for every}  \sigma \in \gal(K/k).
\eqno{(\ref{rem.on.tw.l}.1)}
$$
Conversely, let ${\mathcal F}$ be a  coherent sheaf on $X_K$.
We would like to know whether   ${\mathcal F}\cong F_K$ 
for some coherent sheaf $F$ on $X$. A clear necessary condition is
$$
{\mathcal F}^{\sigma}\cong {\mathcal F}\qtq{for every}  \sigma \in \gal(K/k).
\eqno{(\ref{rem.on.tw.l}.2)}
$$
These conditions are usually not sufficient, and Grothendieck's descent theory
provides the full answer; see \cite[Chap.6]{blr}.

Another approach is to try to understand all  coherent sheaves that
satisfy the naive descent conditions  (\ref{rem.on.tw.l}.2). 
We claim that for  line bundles, this corresponds exactly to the  twisted ones.
(We prove in (\ref{SB.triv.sc.cor})
that we can go from $\bar k/k$ to $k^s/k$ in the definition of
 twisted  line bundles.)

\medskip

{\it Claim \ref{rem.on.tw.l}.3.} Let $X$ be a proper, geometrically reduced and geometrically connected $k$-scheme and ${\mathcal L}$ a line bundle  on $X_{k^s}$.
Then ${\mathcal L}$ is a twisted line bundle  on $X$ iff
$$
{\mathcal L}^{\sigma}\cong {\mathcal L}\qtq{for every}  \sigma \in \gal(k^s/k).
$$

Proof.  Pick  $0\neq F\in T({\mathcal L})$.
For every
$\sigma \in \gal(k^s/k)$ we have
$$
\oplus_i {\mathcal L}\cong F_{k^s}\cong F_{k^s}^{\sigma}\cong \oplus_i {\mathcal L}^{\sigma}.
$$
Thus there is a nonzero map  $ {\mathcal L}\to {\mathcal L}^{\sigma}$.
Applying powers of $\sigma$ until we reach $\sigma^m=1$,
 we get a sequence of nonzero maps
$$
 {\mathcal L}\to {\mathcal L}^{\sigma}\to{\mathcal L}^{\sigma^2}\to  \cdots \to {\mathcal L}^{\sigma^m}\cong {\mathcal L}.
$$
 The composite is an isomorphism, hence all intermediate maps are isomorphisms.

Conversely, let $K/k$ be a finite Galois extension.
Then $ K\otimes_k K$ is a direct sum of
$|\gal(K/k)|$ copies of $K$; projections to the summands correspond to
the bilinear maps  $(a,b)\mapsto a^{\sigma}\cdot b$ where $\sigma\in \gal(K/k)$. 
(See \cite[Secs.III.14-15]{MR0090581}
for a thorough discussion of tensor products of fields.) 

With  $K/k$ as above, let  $p: X_K\to X_k$ be  the natural morphism and $F$ a coherent sheaf on $X_K$. We obtain  that
$$
p^*p_*F=(p_*F) \otimes_k K\cong F\otimes_KK\otimes_k K\cong
\oplus_{\sigma\in \gal(K/k)} F^{\sigma}.
$$
This shows that if  ${\mathcal L}^{\sigma}\cong {\mathcal L}$
for every $\sigma \in \gal(k^s/k)$ then 
 $p_*{\mathcal L}\in T({\mathcal L})$. \qed
\medskip

We can reformulate the above result in terms of the Picard {\em scheme} $\pics(X)$ of $X$. (See  \cite[Lect.19]{mumf66} or \cite[Chap.8]{blr} for  introductions to the Picard  scheme. We will not use any of it in the sequel.)
We get that the class of ${\mathcal L}$ in  $\pics(X)$ is 
invariant under $\gal(k^s/k)$.
Equivalently, $[{\mathcal L}]\in \pics(X)(k)$.
This shows that
$$
\pictw(X)=\pics(X)(k).
\eqno{(\ref{rem.on.tw.l}.4)}
$$
\end{aside}

Next we prove some quite elementary but surprisingly useful results about the category  $T( {\mathcal L})$.

\begin{lem} \label{split.lem} Notation and assumptions are as in
Definition~\ref{T.L.defn}.
 Then all morphisms in $T( {\mathcal L})$ split and there is a unique vector bundle
$E( {\mathcal L})\in T({\mathcal L})$ such that every other member of $T( {\mathcal L})$ is a sum of copies of 
$E( {\mathcal L})$.
\end{lem}

Proof. Let us start over $\bar k$. Since $\Hom({\mathcal L}, {\mathcal L})\cong \bar k$, morphisms $\phi:\oplus_{i=1}^n {\mathcal L}\to \oplus_{j=1}^m {\mathcal L}$ are equivalent to an $m\times n$ matrix
with entries in $\bar k$. Thus kernels and cokernels are also
sums of copies of ${\mathcal L}$ and all morphisms split.

Let now $F_1, F_2\in T({\mathcal L})$ and $\phi:F_1\to F_2$ a morphism.
Then $\ker \phi$ is a coherent sheaf such that 
$(\ker \phi)_{\bar k}=\ker \bigl(\phi_{\bar k}\bigr)$ is a sum of copies of ${\mathcal L}$. Thus
$\ker\phi$ is locally free and it is in $T({\mathcal L})$. Same for
$\coker\phi$. 

An exact sequence of sheaves $0\to F_1\to F_2\to F_3\to 0$  in $T({\mathcal L})$ is classified by an extension class
$\eta\in \ext^1(F_3, F_1)=H^1(X, F_1\otimes F_3^*)$. 
We already know that $\eta_{\bar k}$ is trivial, hence so is $\eta$. 
(Note that we do not claim that if $F_1,F_3\in T({\mathcal L})$ then
$F_2\in T({\mathcal L})$, this is false if $H^1(X, \o_X)\neq 0$.)

Finally let $E({\mathcal L})$ be a vector bundle in $T({\mathcal L})$ with smallest positive rank
and  $F$  any  vector bundle in $T({\mathcal L})$. There is a 
nonzero map $\phi: E({\mathcal L})\to F$ since there is one over $\bar k$. 
Since $\ker \phi\in T({\mathcal L})$, $\phi$ is an injection and
$F\cong E({\mathcal L})+F'$ for some $F'\in T({\mathcal L})$. By induction on the rank we get that 
 $F$ is a sum of copies of 
$E( {\mathcal L})$. \qed

\begin{rem}\label{2.shadows.rem}
 Set $e=\rank E( {\mathcal L})$. Then
$\det E( {\mathcal L})$ is a line bundle on $X$ such that
$$
\bigl(\det E( {\mathcal L})\bigr)_{\bar k}\cong
\det\bigl(\oplus_{1}^e {\mathcal L}\bigr)\cong {\mathcal L}^e.
$$
Thus $\det E( {\mathcal L})$ should be thought of as the
$e$th tensor power of $ {\mathcal L}$ that is defined on $X$.
Thus,   the twisted line bundle   $ {\mathcal L}$
has 2 `shadows' on $X$:
$$
E( {\mathcal L}) \qtq{and} {\mathcal L}^{(e)}:=\det E( {\mathcal L}).
$$
Using these, we will be able to work as if ${\mathcal L} $
were a line bundle on $X$. 
\end{rem}

\begin{cor}  \label{split.lem.cor.1} For every $F\in T({\mathcal L})$,
$\dim_k\End(F)=(\rank F)^2$. Furthermore,
\begin{enumerate}
\item 
$\End\bigl(E( {\mathcal L})\bigr)$ is a skew field,
\item $\End\bigl(\oplus_{i=1}^r E( {\mathcal L})\bigr)$ is an $r\times r$ matrix
algebra over $\End\bigl(E( {\mathcal L})\bigr)$,
\item $F=E( {\mathcal L})$ iff $\End(F)$ is a skew field, 
\item there is a natural correspondence between left ideals  $I\subset \End(F)$ and
direct summands   $G\subset F$, given by
$$
I\mapsto G(I):=\cap_{A\in I}\ker A,\qtq{and}
G\mapsto I(G):=\{A: G\subset \ker A\}.
$$
\end{enumerate}
\end{cor}

Proof. In order to compute  dimensions, note that   
$\End(F)=H^0\bigl(X, \sEnd(F)\bigr)$, where $\sEnd(F) $ denotes the sheaf of endomorphisms, thus
$$
\dim_k \End(F)=h^0\bigl(X_{\bar k}, \sEnd(F)_{\bar k}\bigr)=
h^0\bigl(X_{\bar k}, \sEnd(F_{\bar k})\bigr).
$$
Since $F_{\bar k}$ is a sum of $(\rank F)$ copies of ${\mathcal L}$,
$ \sEnd(F_{\bar k})$ is a sum of $(\rank F)^2$ copies of $\o_{X_{\bar k}}$.
Thus $ \dim_k \End(F)=(\rank F)^2$.

As we noted above,   every endomorphism $\phi:E({\mathcal L})\to E({\mathcal L})$
 is either zero or an isomorphism, thus $\End\bigl(E( {\mathcal L})\bigr)$ is a skew field, and  $\End\bigl(\oplus_{i=1}^r E( {\mathcal L})\bigr)$ is an $r\times r$ matrix
algebra over $\End\bigl(E( {\mathcal L})\bigr)$.
Next  (\ref{split.lem.cor.1}.3) says that 
 $r>1$ iff
there are zero-divisors in the $r\times r$ matrix
algebra.

Finally the natural correspondence between its left ideals and
subspaces of $\End\bigl(E( {\mathcal L})\bigr)^r$ is 
(\ref{split.lem.cor.1}.4).   \qed

\medskip

This is the first appearence of skew fields in our approach to
Severi--Brauer varieties.

\begin{aside} \label{brau.eq.aside} 
In the theory of  central simple algebras,
{\it Brauer equivalence} is obtained by  declaring an algebra $A$ 
to be equivalent to its    matrix algebras $M_{r}(A)$. Thus (\ref{split.lem.cor.1}) says that
the algebras $\End(F)$ are Brauer equivalent to each other for  $F\in T({\mathcal L})$.
\end{aside}

\begin{cor} \label{split.lem.cor.2}
Let $Y$ be another  proper, 
 geometrically reduced and  connected $k$-scheme, and $p:Y\to X$ a morphism. Then $E\bigl(p^* {\mathcal L}\bigr)\cong
p^*E( {\mathcal L})$.
\end{cor}

Proof. It is clear that $p^*E( {\mathcal L})\in T\bigl(p^* {\mathcal L}\bigr)$, the only question is whether $p^*E( {\mathcal L})$ has minimal rank or not.  We have a natural algebra morphism
$\End\bigl(E({\mathcal L})\bigr)\to \End\bigl(p^*E({\mathcal L})\bigr)$.
By (\ref{split.lem.cor.1}) the two have the same dimension, and
$\End\bigl(E({\mathcal L})\bigr)$ is a skew field. Thus 
$\End\bigl(E({\mathcal L})\bigr)\to \End\bigl(p^*E({\mathcal L})\bigr)$ is an isomorphism, and
so $p^*E({\mathcal L})$ is also indecomposable by (\ref{split.lem.cor.1}.3). \qed

\begin{cor} \label{split.lem.cor.3} If $X(k)\neq \emptyset$ then  the
 natural injection
$\pic(X)\into \pictw(X)$ given in (\ref{T.L.defn}.1) is an isomorphism.
\end{cor}

Proof. 
 Pick any twisted line bundle ${\mathcal L}$ on $X$ and 
apply (\ref{split.lem.cor.2}) to $p:\spec k\to X$.
Over $\spec k$ an indecomposable vector bundle has rank 1. Thus
$\rank E({\mathcal L})=1$ and so $E({\mathcal L})\in \pic(X)$ and 
$E({\mathcal L})_{\bar k}\cong {\mathcal L}$ show that
${\mathcal L}\in \pic(X)$.\qed
\medskip

\begin{cor}[Ch\^atelet's theorem] \label{split.lem.cor.4} Let $P$ be a 
Severi--Brauer variety of dimension $n$. Then $P\cong \p^n_k$ iff $P(k)\neq \emptyset$.
\end{cor}

First proof. Assume that $P(k)\neq \emptyset$. We apply (\ref{split.lem.cor.3}) to  $X=P$ and ${\mathcal L}=\o_{P_{\bar k}}(1)$. 
We obtain that $E({\mathcal L})$ is a line bundle  that becomes
isomorphic to $\o_{P_{\bar k}}(1)$ over $\bar k$. Its global sections determine a map $\phi:P\map \p^n_k$ and $\phi$ is an isomorphism over $\bar k$, hence an isomorphism.
The converse  is clear.\qed
\medskip

Second proof. Pick local coordinates $x_1,\dots, x_n$ at a $k$-point of $P$ and consider the restriction map
$$
\begin{array}{rcl}
r_p:H^0\bigl(P, \o_P(-K_P)\bigr)&\to & 
H^0\bigl(P, \o_P(-K_P)\otimes \o_P/(x_1,\dots, x_n)^{n+2}\bigr)\\
&\cong &
\o_P/(x_1,\dots, x_n)^{n+2}.
\end{array}
$$
We know that   $\o_{\p^n}(-K_{\p^n})\cong \o_{\p^n}(n+1)$,
thus $r_p$ is an isomorphism over $\bar k$, hence an isomorphism.
Thus $\o_P(-K_P)$ has a unique section  $s$ such that 
$$
r_p(s)\equiv x_1x_2^n \mod  (x_1,\dots, x_n)^{n+2}.
$$
Over $\bar k$ we can write  $x_i=\bar x_i+ (\mbox{higher terms})$
where the $\bar x_i$ are linear coordinates. Since
$$
 \bar x_1\bar x_2^n \equiv x_1x_2^n \mod  (x_1,\dots, x_n)^{n+2},
$$
we conclude that $(s=0)_{\bar k}= \bigl(\bar x_1\bar x_2^n =0\bigr)$.

For $n\geq 2$ let $H$ denote the closure of the smooth locus of $(s=0)$; it is defined over $k$. 
By the above computation 
$H_{\bar k}$ is the  hyperplane $\bigl(\bar x_1 =0\bigr)$, and hence the linear system $|H|$ maps $P$ isomorphically onto $\p^n$.  \qed

\medskip

\begin{cor}[Wedderburn's theorem]\label{split.lem.cor.5} Over a finite field $\f_q$ every
Severi--Brauer variety is trivial.
\end{cor}

Proof. There are several ways to do this but all have some subtlety.

First proof. It is not hard to prove directly that the number of $\f_q$-points of $P$
is  $1+q+\cdots +q^n$ where $n=\dim P$ and then 
construct an isomorphism $P\cong \p^n$ step-by-step;
see \cite[1.23--26]{ksc} for an elementary approach. 
We could also  use 
(\ref{split.lem.cor.4}) once we have an $\f_q$-point.

Second proof. 
By a theorem of Wedderburn, every finite skew-field is commutative.
Thus  $E(P)$ has rank 1 over $k$ by (\ref{split.lem.cor.1}), so $P$ is trivial. 
The traditional  proof of Wedderburn's theorem is explained in
\cite[Sec.7.2]{MR0171801}, a more powerful argument is in
\cite[Sec.X.7]{MR554237} and \cite[Sec.6.2]{MR2266528}.
\qed

\medskip

\section{First examples  of Severi--Brauer varieties}\label{sec.1.b}

So far the only examples of Severi--Brauer varieties we have are plane conics.
In this section we give 2 methods to construct  new
Severi--Brauer varieties out of a given one, and then 2 constructions of  Severi--Brauer varieties of dimensions 2 and 3.

\begin{say}[Symmetric powers]\label{sym.poerr.defn.1}
Given a variety $X$, let $\sym^n(X)$ denote its $n$th
{\it symmetric power.} That is $\sym^n(X)=X^n/S_n$ is the quotient of the product of $n$ copies on $X$ by the symmetric group $S_n$ permuting the factors.
We can also think of $\sym^n(X)$ as  parametrizing unordered sets of $n$ points
$p_1,\dots, p_n\in X$,  where repetitions are allowed.

One can identify an  unordered set of $n$ points
$p_1,\dots, p_n\in \p^1$ with the unique (up to a multiplicative constant) polynomial of degree $n$ whose roots are $p_1,\dots, p_n$, thus 
 $\sym^n(\p^1)\cong \p^{n}$.
Therefore,  if $C$ is a 1-dimensional
Severi--Brauer variety then $\sym^n(C)$  is an $n$-dimensional
Severi--Brauer variety. 

If $n=2m$ is even then  
$ \sym^{2m}(C) \cong |-mK_C|\cong \p^{2m}$. However,  if $n$ is odd 
and $C(k)=\emptyset$ then $\sym^n(C)$ also has no $k$-points.
Thus 
we get our first examples of  nontrivial  Severi--Brauer varieties of dimension $n>1$.

\smallskip
{\it Aside \ref{sym.poerr.defn.1}.1.} Springer's theorem  says that if a quadric
hypersurface has a point in an odd degree field extension of $k$, then it has a
$k$ point. For 1-dimensional quadrics this is especially easy to prove:
if $P\subset C$ is a point of degree $2r+1$, then any section of
$\o_C(P+rK_C)$ vanishes at a $k$-pont. 
\end{say}

A symmetric power of a variety of dimension $\geq 2$ is always singular, so the symmetric power of a  Severi--Brauer variety of dimension at least 2 can not be a  Severi--Brauer variety. (See, however,  (\ref{kra-sal.thm})  and \cite{k-spsb} for 
a birational study of such symmetric powers.)
Instead we have the following  construction of  larger dimensional
Severi--Brauer varieties from a smaller one. 

\begin{thm} \label{going.up.thm}
Let $P$ be a Severi--Brauer variety, and
$F\in T\bigl(\o_{P_{\bar k}}(1)\bigr)$ a vector bundle. Then
$$
C_F(P)\cong \proj_k \tsum_{m,n=0}^{\infty}  H^0\bigl(P, \sym^{m}(F^*)\otimes 
(\det F)^{\otimes n}\bigr)
\eqno{(\ref{going.up.thm}.1)}
$$
 is a  Severi--Brauer variety  of dimension $\dim P+\rank F$, which  contains $P$
as a subvariety. 
\end{thm}

Note that a summand in (\ref{going.up.thm}.1) is nonzero  iff $m\leq n\cdot \rank F$. 

A conceptual derivation of this formula is given in (\ref{TPSB.min.lem.ef}), which also  proves that $C_F(P)$   is the  unique Severi--Brauer variety 
(up to isomorphism) that  contains $P$ and has  dimension $\dim P+\rank F$.  
\medskip

Proof. Set $r=\rank F$.  We need to check that if $P\cong \p^{\dim P}$ then
$C_F(P)\cong \p^{\dim P+r}$. 
In this case $F\cong \oplus^{r}\o_P(-1)$ and
we can  rewrite the right hand side of
(\ref{going.up.thm}.1) as
$$
 \tsum_{n=0}^{\infty} \tsum_{m=0}^{nr} \sym^{m}(k^{r})\otimes H^0\bigl(P, \o_P(nr-m)\bigr).
\eqno{(\ref{going.up.thm}.2)}
$$
A degree $nr$ homogeneous polynomial in the variables
$x_0,\dots, x_s$ and $y_1,\dots, y_r$ can be uniquely  written as
$$
\tsum_{m=0}^{nr} \tsum_{|I|=m} {\mathbf y}^I f_I({\mathbf x})
\qtq{where} \deg f_I=nr-m.
$$
Thus the sum (\ref{going.up.thm}.2) can be naturally identified with
$$
 \tsum_{n=0}^{\infty} H^0\bigl(\p^{N}, \o_{\p^{N}}(nr)\bigr)
\subset \tsum_{d=0}^{\infty} H^0\bigl(\p^{N}, \o_{\p^{N}}(d)\bigr),
\eqno{(\ref{going.up.thm}.3)}
$$
where $N=\dim P+r$, and two sides of (\ref{going.up.thm}.3)  have the same $\proj$. 

Thinking of $\tsum_n  H^0\bigl(P, (\det F)^{\otimes n}\bigr)$ as a quotient of the sum in (\ref{going.up.thm}.1)  corresponding to $m=0$, 
we see  that
$P\cong \proj_k \tsum_n  H^0\bigl(P,
(\det F)^{\otimes n}\bigr)$ is a subvariety of $C_F(P)$.
\qed

\medskip

If we blow up 6 points in $\p^2$, we get a cubic surface, and the same holds for    Severi--Brauer surfaces. 
One can  reverse this construction by first writing down some 
 cubic surfaces  and then proving that they are birational to  Severi--Brauer surfaces over a field $k$. 
 The key step is to describe the lines 
on a cubic surface, and the  action of the absolute Galois group
$\gal (\bar  k/k)$  on them. 
These tasks are easier for singular cubic surfaces: there are fewer lines and they are easier to locate.

Higher dimensional variants of these examples are discussed in Section~\ref{sec.norm.form}.

\begin{say}[The cubic surface  $S:=(xyz=t^3)$]\label{SB.cubic.geom.say}
The hyperplane section $S\cap (t=0)$ is the union of 3 lines $L_i$ and the surface has $A_2$-singularities at the 3 points where 2 of the lines intersect.
If we resolve these singularities, the birational transforms of the lines $L_i$ and the 6 exceptional curves form a cycle of smooth rational curves.
We can describe it by the dual graph
$$
\begin{array}{ccccccccc}
2 & - & 2 & -\quad 1\quad   - &2 & - &2 \\
| &&&&&& |\\
1 & - & 2 & \rule[1pt]{40pt}{0.4pt} &2 & - &1
\end{array}
\eqno{(\ref{SB.cubic.geom.say}.1)}
$$
where each curve is represented  by minus its self-intersection number, two curves  are disjoint if they are not connected by an edge, and meet at single point if they are  connected by an edge.
The birational transforms of the lines have self-intersection $-1$.  If we contract them we get
the dual graph
$$
\begin{array}{ccccc}
1 & - & 1 & - &1 \\
| &&&& |\\
1 & - & 1 &  - &1
\end{array}
\eqno{(\ref{SB.cubic.geom.say}.2)}
$$
There are 2 ways to choose  3 disjoint $(-1)$-curves from the above 6. Contracting either one of these triples gives $\p^2$. 
\end{say}

We get  interesting variants of these surfaces if we replace
$xyz$ with a cubic form that is irreducible over $k$, but splits into linear factors over $\bar k$. These are the norm forms. 

\begin{defn}[Norm forms]\label{norm.form.defn}
 Let $k$ be a field,  $\alpha\in \bar k$  algebraic  over $k$
and $\alpha_1=\alpha, \alpha_2, \dots, \alpha_n$ its conjugates.
Set $K:=k(\alpha)$. Thus $n=\deg (K/k)$ and
$1, \alpha, \dots, \alpha^{n-1}$ is a $k$-basis of $K$. The {\it norm form} of $\alpha$ (or of $K/k$) is defined as
$$
\norm_{K/k}({\mathbf x})=
\tprod_{i=1}^n\bigl(x_1+\alpha_ix_2+\cdots + \alpha_i^{n-1}x_n\bigr).
\eqno{(\ref{norm.form.defn}.1)}
$$
The product is symmetric in the $\alpha_i$, thus
$$
 \norm_{K/k}({\mathbf x})\in k[x_1,\dots, x_n].
\eqno{(\ref{norm.form.defn}.2)}
$$
The notation suppresses the dependence on the basis 
$1, \alpha, \dots, \alpha^{n-1}$. Using a different basis results in another
norm form that differs by a $k$-linear coordinate change. This will not be important for us.

 Note that $\norm_{K/k}({\mathbf x}) $ is homogeneous of degree $n$,
and has only the trivial zero in $k^n$. Indeed if
$\norm_{K/k}(p_1,\dots, p_n)=0$ then 
$ p_1+\alpha_ip_2+\cdots + \alpha_i^{n-1}p_n=0$ for some $i$, but
the $1, \alpha_i, \dots, \alpha_i^{n-1}$ are linearly independent over $k$.

It is easy to work out that if $\alpha=\sqrt[3]{a}$ and $K=k(\alpha)$, then
$$
 \norm_{K/k}({\mathbf x})=x_1^3+ax_2^3+a^2x_3^3-3ax_1x_2x_3.
\eqno{(\ref{norm.form.defn}.3)}
$$
\end{defn}

\begin{say}[Twisted forms of the cubic surface  $xyz=t^3$]
\label{SB.cubic.arit.say}
Let $k$ be a field and   $K:=k(\alpha)$  a cubic field extension.
 Pick $d\in k$ and consider the cubic surface
$$
 S=S(K, d):=\bigl(\norm_{K/k} (x+\alpha y+\alpha^2z)=dt^3\bigr)\subset \p^3.
 \eqno{(\ref{SB.cubic.arit.say}.1)}
$$
Over $\bar k$ it is isomorphic to $S=(xyz=t^3)$.

It is easy to work out that if $K/k$ is Galois, then the
six $(-1)$-curves in (\ref{SB.cubic.geom.say}.2) form 2 Galois orbits. Contracting either
orbit gives a Severi--Brauer surface.

({\it Aside \ref{SB.cubic.arit.say}.2.} If $K/k$ is not Galois, then the
six $(-1)$-curves  form one Galois orbit, and we have a
degree 6 Del~Pezzo surface of Picard number 1 over $k$.)

As a concrete example, we get the following using (\ref{norm.form.defn}.3).
\smallskip

{\it Example  \ref{SB.cubic.arit.say}.3.} Let $k$ be a field of characteristic $\neq 2,3$ that contains $\sqrt{-3}$. Pick $a,d\in k^{\times}$ 
and assume that  $d$ is not a norm in $k\bigl[\sqrt[3]{a}\bigr]/k$.
Then the cubic surface
$$
(x^3+ay^3+a^2z^3-3axyz=dt^3)
\subset \p^3
$$
is birational to 2  nontrivial  Severi--Brauer surfaces. \qed
\end{say}

Using similar arguments,  Skolem \cite{MR73643}  describes all  Severi--Brauer surfaces.

\begin{thm} \label{SB.cubic.thm}
Let $k$ be a field of  $\chr k\neq 2,3$,
and $K/k$   a cubic extension with discriminant $\Delta_{K/k}$. 
For $\beta\in K^{\times}$ and $d\in k$ consider the cubic surface
$$
T=T(K, \beta, d):=
\bigl(\norm_{K/k}({\mathbf x})+\trace_{K/k}(\beta{\mathbf x})t^2 +d t^3=0\bigr).
\eqno{(\ref{SB.cubic.thm}.1)}
$$
Assume that  $d^2+4\norm_{K/k}(\beta)\neq 0$. Then the following hold. 
\begin{enumerate}\setcounter{enumi}{1}
\item $T$ is rational over $k$ iff  $T(k)\neq \emptyset$.
\item If $T(k)= \emptyset$ and  $\Delta_{K/k}\cdot\bigl(d^2+4\norm_{K/k}(\beta)\bigr)$ is  not a square in $k^{\times}$, then $T$ is birational to a minimal, degree 6 Del~Pezzo surface over $k$.
\item If $T(k)= \emptyset$ and  $\Delta_{K/k}\cdot \bigl(d^2+4\norm_{K/k}(\beta)\bigr)$ is a square in $k^{\times}$, then $T$ is birational to a pair of nontrivial Severi--Brauer surfaces  $P, P^{\vee}$.
\item Every  Severi--Brauer surface over $k$ is obtained this way.  
More precisely, $P$ is obtained from $T(K, \beta, d)$ for some
$\beta, d$ iff $P(K)\neq \emptyset$.
\end{enumerate}
\end{thm}

\begin{aside}
The study of birational maps between Del~Pezzo surfaces was started by Segre \cite{MR0008171}; see \cite{ksc} for an introduction and 
\cite{blanc2024birationalmapsseveribrauersurfaces} for recent applications
to the Cremona groups  $\bir(\p^n)$.
  \end{aside}

An example of a 3-dimensional Severi--Brauer variety is the following.

\begin{prop} Let $Q\subset \p^4$ be a 3-dimensional quadric hypersurface.
Then $\OG(\p^1, Q)$---the variety parametrizing lines in $Q$---is a 3-dimensional  Severi--Brauer variety.
\end{prop}

Note that $\OG(\p^1, Q)(k)=\emptyset$  if $Q(k)=\emptyset$, but there are many examples where $Q(k)\neq\emptyset$, yet $\OG(\p^1, Q)(k)=\emptyset$.
For example, this happens if $k\subset \r$, the $a_i$ are positive and
$Q=\bigl(a_1x_1^2+\cdots+a_4x_4^2=a_0x_0^2)$, since
there are no real lines on an ellipsoid.
\medskip

Sketch of proof. We need to show that $\OG(\p^1, Q)_{\bar k}\cong \p^3_{\bar k}$.
Let $\ell, \ell'$ denote lines on $Q^3$.
We claim that the isomorphism is given by the linear system
$$
\bigl|H_{\ell} :=\{\ell': \ell\cap\ell'\neq \emptyset\}\bigr|.
$$
The key property is to check that $(H_{\ell}^3)=1$. Equivalently,
given 3 disjoint lines $\ell_1,\ell_2,\ell_3\subset Q$, there is a unique line
$\ell'\subset Q$ that meets all 3. To see this note that 
$S=Q\cap \langle \ell_1, \ell_2\rangle$ is a quadric surface
and $p=\ell_3\cap \langle \ell_1, \ell_2\rangle$ is  a point on $S$.
There is thus a unique line $\ell'$ on $S$ that passes through $p$, and
meets  $ \ell_1, \ell_2$. \qed

\section{Twisted linear subvarieties}\label{sec.2}

In this section we  study 
Severi--Brauer varieties using their twisted linear subvarieties.
First we show in (\ref{twisted.basic.thm})   that  twisted linear subvarieties of a Severi--Brauer variety behave very much like linear subspaces of a  vector space.
Then we show another way  to associate a skew field  to a 
Severi--Brauer variety.

\begin{defn} \label{tw.lin.sub.defn}
A subscheme $X\subset P$ of a 
Severi--Brauer variety is called {\it twisted linear} if
$X_{\bar k}$ is a linear subspace of $P_{\bar k}\cong \p^{\dim P}_{\bar k}$.

Let $k^s\subset \bar k$ denote the separable closure. Then 
$P_{k^s}\cong \p^{\dim P}_{k^s}$ by  (\ref{SB.triv.sc.cor}),
hence, if  $X_{\bar k}$ is a linear subspace, then so is
$X_{k^s}$,  since they are both given by the linear forms in the homogeneous ideal of $X_{k^s}$. Thus 
 $X\subset P$ is  twisted linear iff
$X_{k^s}$ is a linear subspace of $P_{k^s}$.

The {\it  twisted linear span} of a closed subscheme  $Z\subset P$
is denoted by  $\langle Z\rangle$. Note that 
$\langle Z\rangle_{\bar k}=\langle Z_{\bar k}\rangle$
which is equivalent to  $\langle Z\rangle_{k^s}=\langle Z_{k^s}\rangle$.
Indeed, $\langle Z\rangle_{k^s}$ is a linear subspace that contains
$Z_{k^s}$, hence $\langle Z\rangle_{k^s}\supset \langle Z_{k^s}\rangle$.
Conversely, we need to show that $ \langle Z_{k^s}\rangle$ is
defined over $k$. It is clearly invariant under $\gal(k^s/k)$, 
thus 
$ \langle Z_{k^s}\rangle$ is defined over $k$ by 
(\ref{weil.lemma}). 

All twisted linear subspaces form a {\it lattice}  
${\mathcal L}(P)$ with intersection and span as the lattice operations.

A map $g:P\map Q$ between Severi--Brauer varieties is called {\it twisted linear} if  it is linear over $\bar k$. That is, the composite
$$
\p_{\bar k}^{\dim P}\cong P_{\bar k}\stackrel{g}{\map}
 Q_{\bar k}\cong \p_{\bar k}^{\dim Q}
\qtq{is linear.}
$$
The key examples of twisted linear maps are projections.
Let $L\subset \p^n$ be a linear subspace. The projection of $\p^n$ 
with center $L$ is the map $\pi_L$ given by  
those sections of $\o_{\p^n}(1)$ that vanish along $L$. 
If $L^{\perp}\subset \p^n$ is a linear subspace complementary to $L$,
then we can think of $\pi_L$ as a map 
$\pi_L:\p^n\map L^{\perp}$. 

Next let  $P$ be a Severi--Brauer variety and
$R\subset P$ a twisted linear subvariety. Although 
$\o_{P}(1)$ does not exist, $\o_{P}(\dim P+1)$ does and
we can define the  {\it projection} 
$\pi_R:P\map R^{\perp}$ of $P$ with {\it center} $R$ as the map 
given by  
those sections of $\o_{P}(\dim P+1)$ that vanish along $R$ with multiplicity $\dim P+1$. It is clear that $R^{\perp}$ is also a Severi--Brauer variety.

We say that two twisted linear subvarieties $Q$ and $ R$
are {\it complementary} if $Q\cap R=\emptyset$ and  
 $\langle Q, R\rangle=P$. Equivalently, if $Q\cap R=\emptyset$ and  
$\dim Q+\dim R=\dim P-1$. 
Then $\pi_R$ restricts to an isomorphism  $Q\cong R^{\perp}$. I particular,
complements of $R$ are all twisted linearly isomorphic to each other.
The projection  $\pi_R:P\map R^{\perp}\cong Q$  has a simple   geometric description: 
 a point $p\in P\setminus R$ is mapped to 
 $\pi_R(p)=Q\cap \langle R, p\rangle$. 
We see in (\ref{TPSB.min.lem.3.cor.2})  that complements exist.
\end{defn}

We have used 
Weil's lemma on the field of definition of a subscheme;
see \cite[I.7.Lem.2]{MR0144898}   or \cite[Sec.3.4]{ksc} for proofs.

\begin{lem}\label{weil.lemma}  Let $X$ be a $k$-scheme, $K/k$ a Galois extension and $Z_K\subset X_K$ a closed subscheme that is invariant under $\gal(K/k)$. Then $Z_K$ is defined over $k$. That is, 
there is a unique closed subscheme $Z_k\subset X$ such that $(Z_k)_K=Z_K$.\qed
\end{lem}

\begin{thm}[Ch\^atelet correspondence]\label{TPSB.cor.2}  
Let $P$ be a Severi--Brauer variety. Then there are natural correspondences between
\begin{enumerate}
\item direct summands  of $F(P)$,
\item left ideals of $\End\bigl(F(P)\bigr)$, and
\item twisted  linear subvarieties of $P$.
\end{enumerate}
 The map between the direct summands and the twisted  linear subvarieties preserves  the lattice structure, that is, intersections, linear spans and complements.   
\end{thm}

Proof.  The equivalence between left ideals  $I\subset \End\bigl(F(P)\bigr)$ and
direct summands   $G\subset F(P)$, given by
$$
I\mapsto G(I):=\cap_{A\in I}\ker A,\qtq{and}
G\mapsto I(G):=\{A: G\subset \ker A\},
\eqno{(\ref{TPSB.cor.2}.4)}
$$
was already discussed in (\ref{split.lem.cor.1}).

To go from direct summands   $G\subset F(P)$ to
twisted  linear subvarieties consider 
$\o_P\to F(P)$ as in (\ref{main.exmp}.4). Set
$$
L(G):=\mbox{zero locus of the section }  \o_P\to F(P)\to F(P)/G.
$$
$\bigl(F(P)/G\bigr)_{\bar k}$ is a sum of copies of $\o_{P_{\bar k}}(1)$, hence
$L(G)$ is a common zero locus of some linear polynomials, thus a
twisted  linear subvariety.
This $G\mapsto L(G)$ is a bijection over $k^s$, hence invertible and  defined over $k$. \qed
\smallskip

Note that the inverse will be clearer once we write the Euler sequence in an invariant way (\ref{TPN.2.say}.4). Restriction gives a surjection
$H^0\bigl(P_{k^s}, \o(1)\bigr)\onto H^0\bigl(L_{k^s}, \o(1)\bigr)$. Its dual is a subspace of $W^*$ in  (\ref{TPN.2.say}.4), hence a
direct summand of $F(P)$.

\begin{cor} \label{TPSB.min.lem.3.cor.2} For every  Severi--Brauer variety $P$ the following hold. 
\begin{enumerate}
\item The minimal twisted linear subvarieties span $P$.
  \item Every twisted linear subvariety $Q\subset P$ has a complement.
\item Let $Q_1\neq  Q_2$  be minimal twisted linear subvarieties. Then  their span 
$\langle Q_1, Q_2\rangle$ contains at least one more minimal twisted linear subvariety
$Q_{12}$. 
\end{enumerate}
\end{cor}

Proof.  We can write $F(P)\cong \oplus_m E_m$ where $E_m\cong E(P)$. Let
$L_m\subset P$ be the  minimal twisted linear subvariety corresponding to $E_i$. Then the $L_m$ span $P$.

Let $F(Q)\subset F(P)$ be the summand corresponding to $Q$. By (\ref{split.lem}) it has complement  $F(P)=F(Q)\oplus F'$, and then $Q':=L(F')$ is a complement of $Q$. 

Let $E_1$ and $E_2$ be the summands that correspond to  $Q_1, Q_2$.
Any isomorphism $\phi:E_1\cong E_2$ gives a summand
$E_\phi$ of $E_1+E_2$, then $L(E_\phi)\subset \langle Q_1, Q_2\rangle$ works.
\qed

\begin{defn}\label{TPSB.cor.3} A Severi--Brauer variety $P$ is 
called {\it minimal} if the  following
equivalent conditions hold.
\begin{enumerate}
\item The only twisted linear subvariety of $P$ is $P$ itself. 
\item $F(P)$ is indecomposable.
\item $\End\bigl(F(P)\bigr)$ is a skew field. 
\end{enumerate}

The equivalence of these follows from (\ref{TPSB.cor.2}.
\end{defn}

\begin{defn-thm} \label{twisted.basic.thm} Let $P$ be a Severi--Brauer variety. Then all minimal (with respect to inclusion) twisted linear subvarieties
are twisted linear isomorphic. Their isomorphism class is denoted by $P^{\rm min}$.

To see the isomorphisms, let   $Q_1, Q_2\subset P$ be 
minimal (with respect to inclusion) twisted linear subvarieties.
Let $R\subset P$ be  a complement to $\langle Q_1, Q_2\rangle$ and
$Q_{12}\subset \langle Q_1, Q_2\rangle$ be a common complement to the $Q_i$ in $\langle Q_1, Q_2\rangle$. Then $\langle Q_{12}, R\rangle$ is a
common complement to the $Q_i$ in $P$. As we noted in (\ref{tw.lin.sub.defn}), both $Q_1,  Q_2$  are isomorphic to  $\langle Q_{12}, R\rangle^\perp$. \qed
\end{defn-thm}

\medskip
{\bf Aside on the Fundamental Theorem of Projective Geometry}
\medskip

\begin{defn} \label{abst.proj.sp.defn}
An abstract  {\it projective space}  ${\mathbf P}:=(P, L)$ 
consists of  a set  $P$ 
(called  points) and  a set  $L$ of subsets of $P$ (called  lines) satisfying the following axioms.
\begin{enumerate}
\item For any two distinct points $p_1,p_2\in P$ there is a unique  line containing them. We denote this line by $\langle p_1, p_2\rangle$.
\item Let $p_1,\dots, p_4$ be 4 distinct points. If
$\langle p_1, p_2\rangle\cap \langle p_3, p_4\rangle\neq \emptyset$
then also $\langle p_1, p_3\rangle\cap \langle p_2, p_4\rangle\neq \emptyset$.
\item  Every line has at least 3 points.
\end{enumerate} 
({\it Comment.} The aim of the definition of  abstract  projective spaces is to axiomatize the properties of linear subspaces of  projective spaces over fields. This leads to a long list of axioms but it turns out that everything can be derived from the above 3.) 

A subset $M\subset P$ is called a {\it linear subspace} if
 $p_1,p_2\in M$ implies that  $\langle p_1, p_2\rangle\subset M$. 
The {\it dimension} of ${\mathbf P}$ is the supremum of the lengths of  chains of linear subspaces
$$
\emptyset\subsetneq M_1\subsetneq \cdots \subsetneq M_d \subsetneq P.
$$
Thus  $\dim {\mathbf P}=0$ if $|P|=1$ and 
$\dim {\mathbf P}= 1$ if $|L|=1$.

Four points $p_1,\dots, p_4$ are {\it coplanar} if 
$\langle p_1, p_2\rangle\cap \langle p_3, p_4\rangle\neq \emptyset$.
Thus axiom (\ref{abst.proj.sp.defn}.2) says that this does not depend on the order of the 4 points.
It is clear that   $\dim {\mathbf P}\leq 2$ iff   any 4 points are coplanar.
\end{defn}

\begin{exmp}\label{subspacelatt.defn}
 Let $F$ be a skew field and $V$ an $F$-vector space.
We denote by   ${\mathcal P}(V)$  the  projective space consisting of the sets
$$
\begin{array}{rcl}
P & := & \{\mbox{1-dimensional subspaces of $V$}\} \qtq{and}\\
L & := & \{\mbox{2-dimensional subspaces of $V$}\}.
\end{array} 
$$
To be precise, a line is the set of all 1-dimensional subspaces
contained in a given 2-dimensional subspace.
It is easy to check that  ${\mathcal P}(V)$  is an abstract  projective space. 
 (Note that this is the dual of the Grothendieck convention that is
common in algebraic geometry.)
\end{exmp}

A fundamental theorem says that 
(\ref{subspacelatt.defn}) gives all  projective spaces of dimension $\geq 3$.
The  theorem is due to  von~Staudt  \cite{vStaudt-2}, but
the first complete proofs are in 
Reye's lectures \cite{zbMATH02723115} (starting with the second edition), using a correction  by Klein \cite{zbMATH02716823}.
Improved versions were given by    Russell \cite{zbMATH02655800},
Whitehead \cite{whitehead},  
 Veblen and Young   \cite{MR1506049}.

\begin{thm}[Fundamental Theorem of Projective Geometry] \label{veb-you.thm} Let  ${\mathbf P}$ be a projective space of dimension $\geq 3$. Then there is a unique  skew field $F$ and an $F$-vector space $V$ such that  ${\mathbf P}\cong {\mathcal P}(V)$. 
\end{thm}

See \cite{MR0172154, MR1153019} for introductions, and
 \cite[Chap.VI]{hodge-ped}, \cite{MR0052795, MR0082463}, or 
 \cite[Chap.5]{MR2264641}
 for full treatments.

Applying this to the lattice ${\mathcal L}(P)$ gives the following. 

\begin{cor} \label{inverse.chat.sec.2}
 Let  $P$ be a Severi--Brauer variety over a field $k$ and 
${\mathcal L}(P)$ the lattice of its twisted linear subspaces.
Assume that $\dim {\mathcal L}(P)\geq 3$.
 Then there is a  unique  skew field $F$ and an $F$-vector space $V$ such that  ${\mathcal L}\cong {\mathcal L}(V)$, the lattice of all linear subspaces of $V$. \qed
\end{cor}

Note that (\ref{inverse.chat.sec.2}) also holds if $\dim {\mathcal L}(P)=2$, but
 there are abstract
 projective planes that do not come from skew fields; see for instance \cite[Chap.6]{MR2264641}.

\section{Splitting fields}\label{sec.split.f}

\begin{defn} Let $P$ be a Severi--Brauer variety of dimension $n$ over $k$.
  A field extension $K/k$ is a {\it splitting field} of $P$ iff
  $P$ has a $K$-point. By (\ref{split.lem.cor.4}) this is equivalent to  $P_K\cong \p^n_K$.
  \end{defn}

Clearly $\bar k$ is a splitting field, but there are much smaller ones.
We start with a weak but useful claim, and then prove the optimal result.

\begin{lem}\label{SB.triv.sc.cor}  Every Severi--Brauer variety has a separable splitting field.
\end{lem}

Proof.  As we discuss in (\ref{stb.ks.say}), every smooth variety
has a separable point. \qed

\begin{thm}\label{TPSB.cor.1} Let $P$ be a  Severi--Brauer variety over an infinite field $k$. Then
  \begin{enumerate}
  \item   The zero-set of a  general section of $T_P$ is a smooth, 0-dimensional subscheme  of degree  $\dim P+1$.
    \item $P$  has a separable splitting field of degree $\dim P^{\rm min}+1$.
    \end{enumerate}
\end{thm}

Proof.  Since $P^{\rm min}\subset P$, we may as well assume that $P$ is minimal.

Fix coordinates $x_0,\dots, x_n$ on $\p^n$. The tangent bundle of $\p^n$
admits a presentation
$$
0\to \o_{\p^n}\stackrel{e}{\to} \oplus_{i=0}^n\o_{\p^n}(1)
\stackrel{\partial}{\to}  T_{\p^n}\to 0,
\eqno{(\ref{TPSB.cor.1}.3)}
$$
Now choose different $c_i\in k$ and
consider the section $(c_0x_0, \dots, c_nx_n)$  of
$\oplus_{i=0}^n\o_{\p^n}(1)$.  Let $\sigma\in H^0(\p^n,  T_{\p^n})$ be its image.
Then   $(s=0)$ is a smooth subscheme consisting of the coordinate vertices.
Thus the zero-set of a  general section of $T_{\p^n}$ is a smooth, 0-dimensional subscheme  $Z$ of degree  $n+1$. If $k$ is infinite, the same holds for
$T_P$.

Let $\emptyset\neq Z'\subset Z$ be a subscheme.
Since $P$ is minimal, the twisted linear span of $Z'$  must equal $P$.
So $\deg Z'\geq n+1$, hence  $Z$ is irreducible (over $k$). 
\qed

\medskip
{\it Aside \ref{TPSB.cor.1}.4.} The algebraic variant of this claim is usually 
 proved using the reduced norm as in \cite[4.5.4]{MR2266528}.

 \medskip

 Another use of $H^0(P, T_P)$ is the following.

\begin{prop} \cite{MR2090670}   \label{kra-sal.thm}
Let $P$ be a Severi--Brauer variety of dimension $n$. 
Then $\sym^{n+1}(P)$ is rational.
\end{prop}

The following proof is taken from  \cite[8]{k-spsb}. Sending a section of $T_P$ to its zero set 
 gives a rational map  $\pi:H^0(P, T_P)\map \sym^{n+1}(P)$.

Let $Z\subset P$ be a reduced 0-cycle of degree $n+1$ whose linear span equals $P$. 
Then  $\pi^{-1}(Z)$ is the linear space
$H^0\bigl(P, T_P(-Z)\bigr)\subset H^0(P, T_P)$ of dimension $n+1$.
Let $V\subset H^0(P, T_P)$ be a general affine-linear subspace of codimension
$n+1$. Then  $\pi|_V:V\map \sym^{n+1}(P)$ is birational. \qed

\begin{aside}[Separable points] \label{stb.ks.say}
Every $k$ variety has points in a finite extension  $k'/k$.
In positive characteristic it is frequently very useful  to have points in 
a finite  and separable extension.

There are at least 3 equivalent versions of the existence of separable points. The first variant is called 
the Separating transcendence basis theorem, proved by MacLane in 1939;
see  \cite[p.18]{MR0144898}, \cite[8.37]{MR571884}, \cite[Sec.VIII.4]{lang-alg} or \cite[p.558]{eis-ca}
  for algebraic proofs.

\medskip
{\it Theorem \ref{stb.ks.say}.1.} Let $K/k$ be a finitely generated field extension of characteristic $p>0$ and of transcendence degree $n$. 
Assume that  $K\otimes_k k^{1/p} $ has no nilpotents (hence, in fact, it is a field).  Then there is a  sub-extension  $K\supset k(x_1,\dots, x_n)\supset k$
such that $K\supset k(x_1,\dots, x_n)$ is finite, separable, and
$ k(x_1,\dots, x_n)\supset k$ is
purely transcendental.

\medskip
{\it Theorem \ref{stb.ks.say}.2.}  
Let $X$ be  a geometrically reduced $k$-variety of dimension $n$.
Then there is a generically finite and separable map $\pi: X\map \p^n_k$. 

(There are sharper variants of this. 
If $X$ is projective, one can choose $\pi: X\to \p^n_k$ to be finite  and separable;
if $X$ is affine one can choose $\pi: X\to \a^n_k$ to be finite and separable.
See \cite[16.18]{eis-ca} for the latter version.)

\medskip
{\it Theorem \ref{stb.ks.say}.3.} A $k$-variety 
 $X$ is geometrically reduced  iff
 $X$ has a smooth point over a finite and separable extension of $k$.
\end{aside}

\section{Non-minimal Severi--Brauer varieties}\label{sec.2.c}

All  Severi--Brauer varieties are built up from
minimal Severi--Brauer varieties in a transparent way.
The key to this is the following result, which shows how to
construct a Severi--Brauer variety from a 
twisted  linear subvariety.

\begin{prop} \label{TPSB.min.lem.3} Let $P$ be a Severi--Brauer variety and
 $Q\subset P$  a  twisted linear subvariety. Let 
$\pi: P\map R:=Q^{\perp}$  be the corresponding projection,  and 
  $\sigma: B_QP\to P$ the blow-up of $Q$. Then the following hold.
\begin{enumerate}
\item The composite $\tilde \pi:= \pi\circ \sigma: B_QP\to  R$ is a
$\p^{\dim Q+1}$-bundle.
\item There is a vector bundle $F\in T(R)$  such that 
$B_QP\cong \p_R\bigl(F^*+\o_R\bigr)$. 
\item $P$ is  uniquely determined by $R$ and $\dim Q$.
  \item  $P$ is birational to $Q\times \p^c$ for $c=\dim P-\dim Q$.
\end{enumerate}
\end{prop}

Proof.   Let $E\subset B_QP$ be the exceptional divisor. 
A geometric fiber of $\tilde \pi$ is a linear space $Q_{\bar k}\subset L\subset P_{\bar k}$ of dimension $=\dim Q+1$, and $E\cap L$ is identified with $Q\subset L$.
 
Pushing forward the injection $\o_{B_QP}\subset \o_{B_QP}(E) $ we get
$$
0\to \o_R\to \tilde \pi_* \o_{B_QP}(E)\to F^*\to 0,
\eqno{(\ref{TPSB.min.lem.3}.5)}
$$
where $F^*$ is a vector bundle of rank $=\dim Q+1$ on $R$.

Working over $\bar k$  shows that
$F^*_{\bar k}\cong \o_R(-1)^{\oplus \dim Q+1}$. Thus
(\ref{TPSB.min.lem.3}.2) holds and (\ref{TPSB.min.lem.3}.5) splits.

Thus $R$ and $\rank F$  determine $B_QP$. There are many ways to see that they also  determine $P$.  We give an explicit formula for $P$ in (\ref{TPSB.min.lem.ef}), but the simplest is to note that  $B_QP$ has Picard number 2, hence at most 2 nontrivial contractions. One of them is
$B_QP\to R$ the other  is $B_QP\to P$. Thus  $B_QP$ uniquely  determines $P$.

Finally (\ref{TPSB.min.lem.3}.4) follows from (\ref{TPSB.min.lem.3}.2). \qed

\begin{say}[Explicit description of $P$] \label{TPSB.min.lem.ef} 
Over $\bar k$ 
$$
\sigma^*\o_{P_{\bar k}}(1)\cong 
\o_{B_QP}(E) \otimes \tilde\pi^*\o_{R_{\bar k}}(-1).
\eqno{(\ref{TPSB.min.lem.ef}.1)}
$$
 Note that $\o_{B_QP}(E)$ is a line bundle, 
but the other 2 are twisted line bundles.
Setting  $c:=\dim P+1$, we get an isomorphism of line bundles over $k$
$$
 \sigma^*\o_{P}(c)\cong \o_{B_QP}(cE)\otimes \tilde\pi^*\o_{R}(c).
\eqno{(\ref{TPSB.min.lem.ef}.2)}
$$
Thus, for every $d\in\n$ we get that
$$
\begin{array}{rcl}
 H^0\bigl(P, \o_P(dc)\bigr)&=&
 H^0\bigl(B_QR,  \o_{B_QP}(dcE)\otimes \tilde\pi^*\o_R(dc)\bigr)\\[1ex]
&=&
 H^0\bigl(R,  \sym^{dc}(F^*+\o_R)\otimes \o_R(dc)\bigr).
\end{array}
\eqno{(\ref{TPSB.min.lem.ef}.3)}
$$
Adding these up for all $d$, the left hand side becomes
$ \tsum_d H^0\bigl(P, \o_P(dc)\bigr)$, whose proj is $P$.
Set $q=\dim Q+1$. 
Adding up the right  hand side for those $d$ such that $q| d$,
and setting $n=dc/q$, 
we get
$$
P\cong \proj_k\Bigl( \tsum_{c|n} H^0\bigl(R, \sym^{nq}(F^*+\o_R)\otimes 
(\det F)^{\otimes n}\bigr)\Bigr).
\eqno{(\ref{TPSB.min.lem.ef}.4)}
$$
The proj does not change if we let the summation run over all $n\geq 0$.
Expanding the symmetric powers now gives that
$$
P\cong \proj_k\Bigl(
\tsum_{n\geq 0}\tsum_{m=0}^{nq}  H^0\bigl(R, \sym^{m}(F^*)\otimes 
(\det F)^{\otimes n}\bigr)\Bigr).
\eqno{(\ref{TPSB.min.lem.ef}.5)}
$$
We already noted in (\ref{going.up.thm}) that this is
equivalent to
$$
P\cong \proj_k\Bigl(
\tsum_{n,m\geq 0} H^0\bigl(R, \sym^{m}(F^*)\otimes 
(\det F)^{\otimes n}\bigr)\Bigr).
\eqno{(\ref{TPSB.min.lem.ef}.6)}
$$
\end{say}

\begin{thm}[Flag transitivity]\label{flag.tarns.thm}
 Let $P_1, P_2$ be two  Severi--Brauer varieties
such that $P_1^{\rm min}\cong P_2^{\rm min}$.
Let 
$$
\emptyset \subsetneq Q^i_1\subsetneq\cdots \subsetneq Q^i_{m_i}\subsetneq P_i
$$
be two maximal chains of twisted linear subvarieties. Assume that
either $\dim P_1\leq \dim P_2$ or 
$m_1\leq m_2$. Then  there is a twisted linear  embedding $\phi: P_1\into P_2$ such that
$\phi(Q^1_j)=Q^2_j$ for every $1\leq j\leq m_1$.
\end{thm}

Proof.
We already know that all minimal twisted linear subvarieties
have the same dimension $i(P)-1$. By (\ref{twisted.basic.thm}),
 $Q^i_{j+1}$ contains a minimal twisted linear subvariety  $R^i_j$ that is disjoint from  $Q^i_j$.  Thus
$$
Q^i_j\subsetneq \langle Q^i_j, R^i_j\rangle \subset Q^i_{j+1}.
$$
By the maximality of the chains $\langle Q^i_j, R^i_j\rangle = Q^i_{j+1}$,
hence $\dim Q^i_{j+1}=\dim Q^i_{j}+ i(P)$. This shows that
$\dim Q^i_{j}=j\cdot i(P)-1$. Thus the 2 assumptions
$(\dim P_1\leq \dim P_2)$ and 
$(m_1\leq m_2)$ are equivalent.
 After replacing $P_2$ by $Q^2_{m_1}$ we may assume that 
$m_1=m_2$; call this common value  $m$. 

Next we use induction on $m$. 
If $m=0$ then $P_1=P_1^{\rm min}\cong P_2^{\rm min}=P_2$, and we are done.
For $m>0$ we may assume that there is a twisted linear  isomorphism $\phi_m: Q^1_m\cong Q^2_m$ such that
$\phi_m(Q^1_j)=Q^2_j$ for every $j$.
By (\ref{TPSB.min.lem.3}), $Q^i_m$ and $\dim P_i$ uniquely determine $P_i$, hence
$\phi_m$ extends to a twisted linear  isomorphism $\phi: P_1\cong P_2$. \qed

\medskip
We are now ready to define Brauer equivalence of
Severi--Brauer varieties. 

\begin{defn-lem}\label{Br.eq.def.lem}
We say that two Severi--Brauer varieties $P_1, P_2$ are 
{\it similar} or {\it Brauer-equivalent,}
denoted by $P_1\sim P_2$, iff the following equivalent conditions hold.
\begin{enumerate}
\item $P_1^{\rm min}\cong P_2^{\rm min}$.
\item The smaller dimensional one  is
isomorphic to a twisted linear subvariety of the other.
\item There is a   twisted linear map  from the larger 
dimensional one onto the other.
\item There is a   twisted linear map  $g:P_1\map P_2$. 
\end{enumerate}
\end{defn-lem}

Proof. We need to check that the four variants are equivalent.

Assume that $\dim P_2\geq \dim P_1$. 
If (\ref{Br.eq.def.lem}.1) holds then (\ref{flag.tarns.thm}) shows that 
$P_1$ is isomorphic to a twisted linear subvariety of $P_2$.
Projecting from a complement of $P_1\subset P_2$ gives a
twisted linear map  from $P_2$ onto $P_1$. Thus
(\ref{Br.eq.def.lem}.1) $\Rightarrow$ (\ref{Br.eq.def.lem}.2) $\Rightarrow$ (\ref{Br.eq.def.lem}.3) while
(\ref{Br.eq.def.lem}.2) $\Rightarrow$ (\ref{Br.eq.def.lem}.4)
and (\ref{Br.eq.def.lem}.3) $\Rightarrow$ (\ref{Br.eq.def.lem}.4) are both clear.

Finally assume (\ref{Br.eq.def.lem}.4).  Let $Q\subset P_1$ denote the locus of indeterminacy of $g$,
and let $P_1^m\subset P_1$ be a minimal  twisted linear subvariety
contained in a complement of $Q$. Then $g(P_1^m)\subset P_2$
is a minimal twisted linear subvariety that is isomorphic to $P_1^m$,
hence (\ref{Br.eq.def.lem}.4) $\Rightarrow$ (\ref{Br.eq.def.lem}.1). \qed

\begin{aside} \label{inter.prod.say}
Let $P$ be a Severi--Brauer variety and
$Q,R\subset P$ complementary twisted linear subvarieties
with normal bundles $N_{Q,P}$ and $N_{R,P}$.
We have  interesting isomorphisms
$$
\p_Q\bigl(N_{Q,P}^*\bigr)\cong Q\times R\cong \p_R\bigl(N_{R,P}^*\bigr).
\eqno{(\ref{inter.prod.say}.1)}
$$
To see this, let $q\in Q$ be a point and $\ell_q\ni q$ a normal line to $Q$.
The  projection of  $\p_Q\bigl(N_{Q,P}^*\bigr) $ to $R$ is given by
$$
\pi_R: (q, \ell_q)\mapsto \langle  Q, \ell_q\rangle \cap R.
\eqno{(\ref{inter.prod.say}.2)}
$$
Thus the coordinate projection  $\pi_Q: Q\times R\to Q$ is a
 non-trivial  $\p^{\dim R}$-bundle, and $\pi_R: Q\times R\to R$ is a
 non-trivial  $\p^{\dim Q}$-bundle.

Note also that  while  $\o_{Q_{\bar k}}(1)$ and $\o_{R_{\bar k}}(-1)$ are twisted  line bundles, 
$$
L\cong q^*\o_{Q_{\bar k}}(1)\otimes r^*\o_{R_{\bar k}}(-1),
\eqno{(\ref{inter.prod.say}.3)}
$$
is an actual line bundle. A similar formula leads to the definition of the dual 
Severi--Brauer variety in (\ref{dial.defn.1}).
\end{aside}

\section{Twisted linear systems}\label{sec.6}

If $L$ is a line bundle on $X$, then $|L|:=\p\bigl(H^0(X, L)^{\vee}\bigr)$ is the
{\it linear system} associated to $L$. 
If ${\mathcal L}$ is a twisted line bundle on $X$, then 
there does not seem to be any way of defining $H^0(X, {\mathcal L}) $.
However, one can define the 
{\it twisted linear system} $|{\mathcal L}|$ associated to ${\mathcal L}$; it is a Severi--Brauer variety. We discuss three ways of constructing
 $|{\mathcal L}|$. The first---most general one---is equivalent to 
viewing a Severi--Brauer variety as a Galois cohomology class.
The second relies on the existence of the Hilbert scheme. 
The third approach is quite elementary but it does not give the most general result.

\begin{aside}[First construction]\label{|L|.const.1}
  Let $X$ be a proper  $k$-scheme, and   ${\mathcal G}$ a
simple, coherent sheaf 
on $X_{k^s}$ such that ${\mathcal G}^{\sigma}\cong {\mathcal G}$ for every $\sigma\in \gal(k^s/k)$.
As we noted in Aside~\ref{rem.on.tw.l}, this holds
if     ${\mathcal G}={\mathcal L}$ is a twisted line bundle on $X$.

Since  ${\mathcal G}$ is
simple, the isomorphism ${\mathcal G}^{\sigma}\cong {\mathcal G}$
is unique up-to a multiplicative scalar, so we get well-defined
isomorphisms  
$$
j_{\sigma}:\p\bigl(H^0(X_{k^s}, {\mathcal G}^{\sigma})^{\vee}\bigr)\cong
\p\bigl(H^0(X_{k^s}, {\mathcal G})^{\vee}\bigr).
$$
These define a $k$-scheme structure on 
$\p\bigl(H^0(X_{k^s}, {\mathcal G})^{\vee}\bigr) $.
\end{aside}

\begin{say}[Second construction]\label{|L|.const.2}
 Let $X$ be a proper, geometrically connected  and normal  $k$-scheme. Let  ${\mathcal L}$ be a line bundle 
on $X_{k^s}$ such that ${\mathcal L}^{\sigma}\cong {\mathcal L}$ for every $\sigma\in \gal(k^s/k)$.

Let $|{\mathcal L}|$ denote the
irreducible component of the Hilbert scheme  of $X$
parametrizing subschemes $H\subset X$ such that $H_{k^s}$ is in the linear system $|{\mathcal L}|$.

By our assumption   $|{\mathcal L}|$ is invariant under the
Galois group $\gal(k^s/k)$, so it is naturally a $k$-variety
by (\ref{weil.lemma}), hence 
 a Severi--Brauer variety.
\end{say}

\begin{say}[Third  construction using Veronese embeddings] \label{verones.say}
Let $Y$ be a 
  $k^s$-scheme, and $M$ a line bundle on $Y$.
Fix a natural number $e>0$. Sending $s\in H^0(Y, M)$ to
$s^e\in H^0(Y, M^e)$ descends to an embedding
$$
v_{e,M}: |M|\into |M^e|.
\eqno{(\ref{verones.say}.1)}
$$
We apply this to  $Y:=X_{k^s}$, 
$M:={\mathcal L}$,    and $e>0$ such that
there is a line bundle $L^{(e)}$ on $X$ satisfying
$(L^{(e)})_{k^s}\cong {\mathcal L}^e$. 
(We can take  $e=\rank E({\mathcal L})$ by Remark~\ref{2.shadows.rem}.)  Then
$$
v_{e,{\mathcal L}}: |{\mathcal L}|\into |{\mathcal L}^e|\cong |L^{(e)}|_{k^s}.
\eqno{(\ref{verones.say}.2)}
$$
We claim  that  the image of $v_{e,{\mathcal L}} $---which, by construction, is defined by equations over $k^s$---can be defined by  equations over $ k$.
If this holds then the resulting $k$-variety is  $|{\mathcal L}| $. 

We leave the general setting to the reader, but 
discuss   $\bigl(P, {\mathcal L}=\o_{P_{k^s}}(1)\bigr)$; 
 the case that we need later.
We rely on the following simple observation:

A hypersurface $H\subset \p^n$ of degree $e$ is an $e$-fold hyperplane iff every point of it has multiplicity $e$.

Set $n:=\dim P$. We now that  $\o_P(n+1)\cong {\mathcal L}^{(n+1)}$ is a line bundle on $P$, and consider the universal hypersurface with projection 
$$
{\mathbf H}\subset |\o_P(n+1)|\times P \stackrel{\pi}{\to} |\o_P(n+1)|.
\eqno{(\ref{verones.say}.3)}
$$
It is easy to see that  $(H, x)\mapsto \mult_xH$ is an upper semi-continuous function on $ {\mathbf H}$. Thus, for every $c$  there is a closed $k$-subscheme
$W_c\subset |\o_P(n+1)|$ parametrizing those hypersurfaces $H$ of degree $n+1$
such that $\mult_xH\geq c$ for every $x\in H$. This completes the construction in  (\ref{verones.say}), since
$$
v_{n+1,{\mathcal L}}\bigl(|{\mathcal L}|\bigr)=(W_{n+1})_{k^s}.
\eqno{(\ref{verones.say}.4)}
$$
Thus $|{\mathcal L}|:=W_{n+1}$ is the sought after $k$-structure in this case.

Note that this method works for
 $(X, {\mathcal L})$ if  every member of $|{\mathcal L}|$ is generically reduced and $\pic(X_{k^s})$ is torsion free. 
\end{say}

This special case already allows us to define duals.

\begin{defn}\label{dial.defn.1}
The {\it dual}  of  a Severi--Brauer variety $P$ is  defined as $P^{\vee}:=|\o_P(1)|$.
Thus $P^{\vee}$ is a Severi--Brauer variety  that
parametrizes all hyperplanes in $P$. 

There is a natural isomorphism 
$$
\iota_{\mathcal P}:P\cong (P^{\vee})^{\vee}=|\o_{P^{\vee}}(1) |^{\vee}\qtq{given by} x\mapsto \{H: H\ni x\}.
\eqno{(\ref{dial.defn.1}.1)}
$$ 
Thus $P\mapsto P^{\vee}$ is a duality for Severi--Brauer varieties.

It is important to note that the incidence divisor
  $\{(x,H): x\in H\}$ on $P\times P^{\vee} $ is defined over $k$, hence 
$\o_{P\times P^{\vee}}(1,1)$ is a (non--twisted) line bundle on $P\times P^{\vee} $.

\end{defn}

\begin{rem}  \label{ten.M.and.|L|}
Let $M$ be a line bundle on $X$ and
$s:\o_X\to M$ a  section that does not vanish on any irreducible component of $X$. Then tensoring with $s$ gives 
$$
E(s):E({\mathcal L}) \to E({\mathcal L}\otimes M)\qtq{and}
s^{(e)}: L^{(e)}\to (L\otimes M)^{(e)}.
$$
Thus $s^{(e)}$ descends to a 
twisted linear map
$$
\phi_s:  | {\mathcal L}|\into |{\mathcal L}\otimes M|.
$$
In particular,  $| {\mathcal L}|\sim |{\mathcal L}\otimes M|$ by
(\ref{Br.eq.def.lem}).
\end{rem}

\begin{defn}\label{tw.lb.defn} Assume that $X$ satisfies one of the assumptions in 
(\ref{|L|.const.1}--\ref{verones.say})
 and let ${\mathcal L}$ be a
 twisted line bundle on $X$ that  is   generated by global sections.

As in (\ref{dial.defn.1}.1)  there is a natural morphism 
$$
\iota_{\mathcal L}:X\to |{\mathcal L}|^{\vee}\qtq{given by} x\mapsto \{H: H\ni x\}.
\eqno{(\ref{tw.lb.defn}.1)}
$$ 
Note that $\iota_{\mathcal L}^*\o_{|{\mathcal L}|^{\vee}}(1)\cong {\mathcal L}$, hence, by (\ref{split.lem.cor.2}),
$$
E({\mathcal L})=\iota_{\mathcal L}^*E\bigl(|{\mathcal L}|^{\vee}\bigr).
\eqno{(\ref{tw.lb.defn}.2)}
$$
If  ${\mathcal L}$ is not
 base point free, but $|{\mathcal L}|\neq\emptyset$ and $X$ is smooth,  one can still define a rational map 
$$
\iota_{\mathcal L}:X\map |{\mathcal L}|^{\vee}\qtq{given by} x\mapsto \{H: H\ni x\}.
\eqno{(\ref{tw.lb.defn}.3)}
$$ 
Note that   $\iota_{\mathcal L}:X\map |{\mathcal L}|^{\vee} $ does not depend on the fixed part of
the linear system $|{\mathcal L}|$. However, as long as 
$|{\mathcal L}|$ is {\it mobile,} that is, 
the base locus of $|{\mathcal L}|_{k^s}$ 
has codimension $\geq 2$, the isomorphism
(\ref{tw.lb.defn}.2)  holds.

 More generally, let $\phi:X\map Y$ be a  map between proper, 
geometrically connected  and geometrically reduced  
 varieties and ${\mathcal L}_Y$ a 
twisted line bundle on $Y$. Assume that either $\phi$ is a morphism or
$X$ is smooth. Then $\phi^*{\mathcal L}_Y$ is a twisted line bundle on $X$, and
the induced pull-back map  $\phi^*: |{\mathcal L}_Y|\map |\phi^*{\mathcal L}_Y|$ is twisted linear. Thus, by
(\ref{Br.eq.def.lem}.4),
$$
|\phi^*{\mathcal L}_Y|\sim |{\mathcal L}_Y|.
\eqno{(\ref{tw.lb.defn}.4)}
$$
Working  with rank 1 reflexive sheaves  as in
Definition~\ref{vardefs.1} shows that 
(\ref{tw.lb.defn}.4)  holds provided $X,Y$ are normal and $ \phi^*{\mathcal L}_Y$ makes sense, that is,
 $\phi(X)$ is not contained in the base locus of $|{\mathcal L}_Y|$. 
\end{defn}

Next we come to the definition of products. We start with the
standard version and then translate it into other forms.

\begin{defn}[Products] \label{prod.defn}
Here are 5 versions of the definition of products.
The first 2 variants recall known definitions from geometry and algebra.

\medskip
(Polarized pair version \ref{prod.defn}.1)  Let $X, Y$ be proper
$k$-varieties and ${\mathcal L}_X, {\mathcal L}_Y$ twisted line bundles on them.
Their product is
$$
\bigl(X\times Y, {\mathcal L}_{X\times Y}:=\pi_X^*{\mathcal L}_X\otimes \pi_Y^*{\mathcal L}_Y\bigr),
$$
where $\pi_X, \pi_Y$ are the coordinate projections.

\medskip
(Algebra version \ref{prod.defn}.2)  It is clear that
$\pi_X^*E({\mathcal L}_X)\otimes \pi_Y^*E({\mathcal L}_Y)\in T({\mathcal L}_{X\times Y})$. It is not necessarily
of minimal rank in $T({\mathcal L}_{X\times Y})$, but we have that
$$
\End\bigl(\pi_X^*E({\mathcal L}_X)\otimes \pi_Y^*E({\mathcal L}_Y)\bigr)
\cong \End\bigl(E({\mathcal L}_X)\bigr)\otimes\End\bigl(E({\mathcal L}_Y)\bigr).
$$

We can now turn (\ref{prod.defn}.1) into a definition of products of linear systems.

\medskip
(Linear system  version \ref{prod.defn}.3) Following the previous examples, we  set
$$
|{\mathcal L}_X|\otimes |{\mathcal L}_Y|:=|{\mathcal L}_{X\times Y}|=\bigl| \pi_X^*{\mathcal L}_X\otimes \pi_Y^*{\mathcal L}_Y\bigr|.
$$  

Thinking of a Severi--Brauer variety as a twisted linear system
gives the following.

\medskip
(Severi--Brauer version \ref{prod.defn}.4) A Severi--Brauer variety $P$
is isomorphic to  $|\o_{P^{\vee}}(1)|$. This suggests that the product  
of  Severi--Brauer varieties $P,Q$  should be  defined as
$$
 P\otimes Q:=|\o_{P^{\vee}\times Q^{\vee}}(1,1)|, \qtq{which is isomorphic to}
|\o_{P\times Q}(1,1)|^{\vee},
$$
where $ \o_{P\times Q}(a,b):=\pi_P^*\o_P(a)\otimes\pi_Q^*\o_Q(a) $
using the coordinate projections of $P\times Q$. 
We will write 
$P^{\otimes m}:=P\otimes \cdots \otimes P$  ($m$ factors).
\medskip

One can reformulate (\ref{prod.defn}.4) as follows.

\medskip
(Segre embedding version \ref{prod.defn}.5) Let $P_i$ and $Q$ be 
Severi--Brauer varieties.  Then $Q\cong P_1\otimes\cdots\otimes P_r$ iff there is an embedding
$$
j:P_1\times\cdots\times P_r\into Q
$$
that becomes the Segre embedding over $\bar k$. 
\medskip

To see that (\ref{prod.defn}.4) and (\ref{prod.defn}.5) are equivalent,
note that $j$ gives a twisted linear pull-back map
$j^*: Q^{\vee}\to |j^*\o_Q(1)|\cong P_1^{\vee}\times\cdots\times P_r^{\vee}$.
This is an isomorphism over $\bar k$ hence also an isomorphism over $k$.
\end{defn}

\begin{lem} \label{powers.of.L}
 Let $X$ be a proper, geometrically reduced and connected   $k$-variety,
and ${\mathcal L}_i$  base-point free, twisted line bundles on $X$. 
Then $|\otimes_i{\mathcal L}_i|\sim \otimes_i|{\mathcal L}_i|$,
that is, the linear system associated to the  tensor product 
$\otimes_i{\mathcal L}_i $
is similar to the product of the linear systems $|{\mathcal L}_i| $. 
In particular, $|{\mathcal L}^m|\sim |{\mathcal L}|^{\otimes m} $.
\end{lem}

Note that we do not claim isomorphism, only similarity,
as in (\ref{Br.eq.def.lem}).
\medskip

Proof.  Consider the diagonal embedding  $\delta: X\into X^m$,
and the twisted line bundle  ${\mathcal M}:=\otimes_i \pi_i^*{\mathcal L}_i$, where the $\pi_i$ are the coordinate projections.  
Then $\delta^*{\mathcal M}=\otimes_i{\mathcal L}_i$ and, using first
(\ref{tw.lb.defn}.4) and then (\ref{prod.defn}.4), we get that 
$|\otimes_i{\mathcal L}_i|=|\delta^*{\mathcal M}|\sim |{\mathcal M}|\sim 
 \otimes_i|{\mathcal L}_i|$. \qed

\begin{cor} \label{powers.of.L.cor} Let  $p:P\map Q$ be a rational map between 
Severi--Brauer varieties. Then $p^*\o_Q(1)\cong \o_P(m)$ for some  $m$, and 
$Q\sim  P^{\otimes m}$.
\end{cor}

Proof. 
By (\ref{tw.lb.defn}.4)  we get
a linear embedding $|\o_Q(1)|\into |\o_P(m)|$. The latter is
similar to $|\o_P(1)|^{\otimes m}$ by (\ref{powers.of.L}). Thus
$|\o_Q(1)|\sim |\o_P(m)|\sim |\o_P(1)|^{\otimes m}$ and taking duals gives that
$Q\sim P^{\otimes m}$.\qed

\medskip

\begin{exmp}[Cremona transformations I]\label{crem.tr.exmp}
  The simplest birational transformation of $\p^n$ is
  $$
  (x_0{:}\cdots{:} x_n)\mapsto (x_0^{-1}{:}\cdots{:} x_n^{-1}).
  \eqno{(\ref{crem.tr.exmp}.1)}
  $$
  After multiplying by $x_0\cdots x_n$, it is given by the linear system of degree $n$ polynomials that vanish with multiplicity $n-1$ at the coordinate vertices.

  Thus if $P$ is a Severi--Brauer variety of dimension $n$, and
  $Z\subset P$ a 0-dimensional, spanning subscheme of degree $n+1$, then
  $\bigl|\o_P(n)(-(n-1)Z)\bigr|$ defines a birational map
  $P\simb P^\vee$. 
\end{exmp}

\begin{exmp}[Cremona transformations II]\label{crem.tr.exmp.2}
  A more precise version is the following.

  Let $K/k$ be an extension of degree $n+1$ such that $P_K\cong \p^n_K$. 
  Let $H_0\subset P_K$  and $G_0\subset P_K^\vee$ be  hyperplanes whose conjugates  $H_i$  (resp.\ $G_i$) are linearly independent. Set   $H:=\sum H_i$  and $G:=\sum G_i$.

  Let $\pi,\pi^\vee$ denote the coordinate projections of $P\times P^\vee$. By (\ref{dial.defn.1})  $\o_{P\times P^\vee}(1,1)$ is a line bundle,
  and the $n+1$ conjugates of $\pi^*H_0+(\pi^\vee)^*G_0$ span a linear system that is defined over $k$. The intersection of all these conjugates  $\Gamma_{H,G}$
  is the graph of a Cremona transformation  $\phi_{H,G}:P\map P^\vee$ such that
  $\phi_{H,G}^*(G_0)=H-H_0$.

  Indeed, this claim can be checked over $K$. We can then choose
  $H_i=(x_i=0)$, $G_i=(y_i=0)$, and then $\Gamma_{H,G}$
  is given by the equations
  $x_0y_0=\cdots=x_ny_n$.  This is the closed graph of (\ref{crem.tr.exmp}.1).

  Let now $H'_0\subset P_K$ be another hyperplane. Then
  $$
  (\phi_{H,G}\circ \phi_{G,H'})^*(H'_0)=\phi_{H,G}^*(G-G_0)=(n-1)H+H_0.
  $$
  After canceling the common $(n-1)H$, we get that
  $\phi_{H,G}\circ \phi_{G,H'}$ is an automorphism of $P$ that maps
  $H_0$ to $H'_0$.
  \end{exmp}

The definition of $ |{\mathcal L}|$ can be extended to coherent sheaves.

\begin{defn}\label{|F|.defn} Let $F$ be a coherent sheaf on a Severi--Brauer variety $P$, and
  ${\mathcal L}$ a twisted line bundle. Then one can define a
  Severi--Brauer variety
  $$
  |F\otimes \mathcal L|\qtq{such that}
  |F\otimes \mathcal L|_{k^s}=\p\bigl(H^0(P_{k^s}, F_{k^s}\otimes \mathcal L)^\vee\bigr).
  $$
  To do this, let  $A$ be a very ample line bundle on $P$ such that
  $F$ is a quotient  $\oplus_i A^{-1}\onto F$. If $r\gg 1$, then
  $$
  H^0(P_{k^s}, \oplus_i A^{r-1}\otimes \mathcal L)\to H^0(P_{k^s}, F\otimes A^r\otimes\mathcal L)
  $$
  is surjective.
  Since  $ \oplus_i A^{r-1}\otimes \mathcal L$ is a direct sum of twisted line bundles, we get a Severi--Brauer structure on the
  projectivization of $H^0(P_{k^s}, \oplus_i A^{r-1}\otimes \mathcal L)$, which in turn puts a 
Severi--Brauer structure on
the projectivization of $H^0(P_{k^s}, F\otimes A^r\otimes\mathcal L) $.

Choosing a general  $\o_P\into A^r$ then realizes
$H^0(P_{k^s}, F\otimes\mathcal L)$ as a
linear subspace of  $H^0(P_{k^s}, F\otimes A^r\otimes\mathcal L)$.
\end{defn}

\section{The Brauer group}\label{sec.7}

\begin{thm}  The product defined in (\ref{prod.defn}.4) makes the
Brauer equivalence classes of Severi--Brauer varieties into a commutative group. The inverse of $P$ is $P^{\vee}$.
\end{thm}

Proof. The  product defines a commutative monoid structure on
the isomorphism classes of Severi--Brauer varieties.
If  $P_1\map P_2$ and $Q_1\map Q_2$ are twisted linear maps, then so is
$P_1\otimes Q_1\map P_2\otimes Q_2$, hence  we get a
commutative monoid on the Brauer equivalence classes.

As we noted in (\ref{dial.defn.1}), 
$\o_{P\times P^{\vee}}(1,1)$ is a (non--twisted) line bundle, hence
$|\o_{P\times P^{\vee}}(1,1)|$  is a trivial Severi--Brauer variety.
This says that 
 $P\otimes P^{\vee}\sim \p^0$. \qed

\begin{defn} \label{br.gp.defn}
The group defined above is  called the {\it Brauer group} of the field $k$, and denoted by $\br(k)$.   It was originally defined in terms of
tensor products of (Brauer equivalence classes of)
central simple $k$-algebras; see (\ref{brau.eq.aside}). Nowadays it is most commonly defined as
$H^2(k, \gm)$, the second Galois cohomology group of $\gm$. 
See  \cite{MR2266528} for these aspects.

The order of $P$ in
$\br(k)$ is traditionally called that {\it period} of $P$.
Since $P\sim P^{\rm min}$, they have the same period.
\end{defn}

\begin{prop}\label{period.prop} The period of $P$ equals the smallest  $m>0$
such that  $\o_P(m)$ is a (non-twisted) line bundle on $P$.
Thus  the period of $P$ divides  $\dim P+1$ and $\br(k)$ is a torsion group. 
\end{prop}

Proof. By (\ref{powers.of.L}), $|\o_P(m)|\sim |\o_P(1)|^{\otimes m}
\sim(P^{\vee})^{\otimes m}\sim(P^{\otimes m})^{\vee}$. 
As we noted in Example~\ref{main.exmp}, $\o_P(\dim P+1)\cong \o_P(-K_P)$
is a (non-twisted) line bundle on $P$. 
\qed

\begin{prop}[Amitsur's theorem]\label{ami.cor}
 Let $P, Q$ be Severi--Brauer varieties. 
The following are equivalent.
\begin{enumerate}
\item $Q$ is similar to
$P^{\otimes m}$ for some $m$.
\item $P\times \p^{\dim Q}\simb P\times  Q$.
\item There is a rational map  $p:P\map Q$.
\end{enumerate}
\end{prop}

Proof. Assume (\ref{ami.cor}.1) and let $K$ be the function field of $P$. Then $P_K$ has a $K$-point
(the generic point of $P$) thus $P_K$ is trivial and so is
$Q_K\sim P_K^{\otimes m}$. 
Thus $Q_K$  is isomorphic to $\p^{\dim Q}_K$ by (\ref{split.lem.cor.4}), hence
(\ref{ami.cor}.2) holds.  The implication (\ref{ami.cor}.2) $\Rightarrow$ (\ref{ami.cor}.3) is clear.
Finally assume that there is a rational map  $p:P\map Q$. Then 
$p^*\o_Q(1)\cong \o_P(m)$ for some $m$ hence
$Q\sim  P^{\otimes m}$ by (\ref{powers.of.L.cor}). \qed

\medskip

Frequently (\ref{ami.cor}.1) $\Leftrightarrow$  (\ref{ami.cor}.2) is stated in the following equivalent form.

\begin{cor} The kernel of the base-change map
$\br(k)\to \br\bigl(k(P)\bigr)$ is the subgroup generated by $P$. \qed
\end{cor} 

\begin{rem} \label{ami.cor.cor} Applying  (\ref{ami.cor}.1) $\Rightarrow$ (\ref{ami.cor}.2)
twice we obtain that if 
 $P, Q$ generate the same subgroup of $\br(k)$ then
$$
P\times \p^{\dim Q}\simb P\times  Q\simb \p^{\dim P}\times Q.
$$
That is, $P$ and $Q$ are stably birational to each other.
As a special case,
$$P\times  P\simb  P\times \p^{\dim P}.$$

 It is an unsolved problem whether  the stronger assertion
$P^{\rm min}\simb Q^{\rm min}$
holds.

For Severi--Brauer varieties that  become trivial after a cyclic Galois extension, the birational isomorphism
$P^{\rm min}\simb Q^{\rm min}$ follows from (\ref{normform.SB.thm.compl}.2).

Besides Cremona transformations as in (\ref{crem.tr.exmp})
the only other known general birational isomorphisms among  Severi--Brauer varieties  is the following.
\end{rem}

\begin{exmp} \label{tregub.map} \cite{tregub}
Let $V$  be an odd dimensional $k$-vectorspace and $\chr k\neq 2$.
We can think of $\wedge^2V$ is skew-symmetric matrices and a general
skew-symmetric matrix has a unique eigenvector with eigenvalue 0.
This gives a rational map  $\wedge^2 V\map \p(V)^{\vee}$.

For an  even dimensional 
Severi--Brauer variety $P$, define
$\wedge^2 P\subset |\o_{P^\vee\times P^\vee}(1,1)|$ as the sections that vanish on the diagonal. Thus there is a rational map
$\psi: \wedge^2 P\map P$.
This map is not twisted linear, but its fibers are linear
subspaces. Thus if $Q\subset \wedge^2 P$ is a  minimal  twisted linear
subspace,  then $\psi|_Q:Q\map P$ is birational.
Therefore
$(P^{\otimes 2})^{\rm min}\cong Q \simb P^{\rm min}$.

There are many numbers $n$ such that $2$ and $-1$ generate the
multiplicative group  $(\z/n)^{\times}$. Thus we obtain that if
$n$ is such a number, and 
 $P, Q$ generate the same subgroup of order $n$ of $\br(k)$, then
$P\simb Q$. 
\end{exmp}

\section{Index of a Severi--Brauer variety}\label{sec.8}

\begin{defn}[Index of a  variety] Let $X$ be a proper  $k$-scheme.
The {\it index} of $X$ is the gcd of the degrees of all 
 0-cycles on $X$.
It is denoted by  $\ind(X)$.  

For smooth, proper varieties the index is a birational invariant by (\ref{nishi.lem}).
More general indices are defined in \cite{k-elw}, building on  \cite {MR3383601}. 
\end{defn}

A short proof of the following is in \cite[p.183]{ksc}.

\begin{lem}[Nishimura lemma]\label{nishi.lem} Let $X\map Y$ be a rational map
between $k$-varieties. Assume that $X$ has a smooth $k$-point and $Y$ is proper. Then $Y$ has a $k$-point. \qed
\end{lem}

\begin{lem}\label{submult.index.lem}
 Let $X_k$ be a normal, proper  $k$-variety
and $K/k$ a finite field extension. Then
$\ind(X_k)\mid \deg (K/k)\cdot \ind(X_K)$.
\end{lem}

Proof. $p:X_K\to X_k$ is a finite map whose degree equals
$\deg (K/k) $. Thus  if $Z$ is a 0-cycle on $X_K$, then
$\deg (p_*Z)=\deg (K/k)\cdot\deg Z$. \qed

\begin{lem} \label{ind.of.prods}
Let $P,Q$ be Severi--Brauer varieties. Then
$\ind(P\otimes Q)$ divides $\ind(P)\cdot \ind(Q)$
and  $\ind(P^{\otimes m})$ divides  $\ind(P)$.
\end{lem}

Proof. Let $Z_P\subset P$ and $Z_Q\subset Q$ be 0-cycles.
Then $Z_P\times Z_Q\subset P\times Q$ is a 0-cycle
whose degree is $\deg(Z_P)\cdot\deg(Z_Q)$. By
(\ref{prod.defn}.5) it is also a 0-cycle on $P\otimes Q $.

The diagonal embedding shows that $Z_P$ is also a 0-cycle on
$P^m\subset P^{\otimes m}$. \qed

\begin{cor} If $P, Q$ generate the same subgroup of $\br(k)$ then they have the same index. \qed
\end{cor}

\begin{thm} \label{SB.index.thm}
Let $P$ be  a Severi--Brauer variety. Then
\begin{enumerate}
\item $\ind(P) \mid \dim P+1$, 
\item  $\ind(P)=\ind( P^{\rm min})$,
\item  $\ind(P)=\dim P^{\rm min}+1$, and
\item  $P$ contains  smooth, 0-dimensional subschemes  whose  degree equals  $\ind(P)$.
\end{enumerate} 
\end{thm}

Proof.  By (\ref{TPSB.cor.1}), $P$ contains a 0-cycle of degree $\dim P+1$, so (\ref{SB.index.thm}.1) holds by definition. 
By (\ref{TPSB.min.lem.3}.4) $P\simb P^{\rm min}\times \p^c$
hence any 0-cycle on $P$ (resp.\ $P^{\rm min}$) yields a
0-cycle of the same degree on $P^{\rm min}$ (resp.\ $P$).
This shows (\ref{SB.index.thm}.2).

 By (\ref{TPSB.cor.1}), $P^{\rm min}$, and hence $P$,
contain smooth, $0$-dimensional subschemes  of degree $\dim P^{\rm min}+1$. Thus
$\ind(P)\leq \dim P^{\rm min}+1$.

Finally let $Z\subset P$ be a 0-dimensional reduced subscheme.
By (\ref{|F|.defn}), $|\o_Z(-1)|$ is a Severi--Brauer variety similar to $P$.
Its dimension is $\deg Z-1$, so $\deg Z$ 
 is divisible by
$\dim P^{\rm min}+1$.  This  completes (\ref{SB.index.thm}.3) and proves (\ref{SB.index.thm}.4).  \qed

\begin{cor}\label{index.drops.pow.cor} Let $P$ be a Severi--Brauer variety and $m$ a divisor of $\ind(P)$. Then
$\ind(P^{\otimes m})$ divides  $\tfrac1{m}\ind(P)$.
\end{cor}

Proof. Write $m=m_1m_2$. If $\ind(P^{\otimes m})\nmid \tfrac1{m}\ind(P)$
then either $\ind(P^{\otimes m})\nmid \tfrac1{m_2}\ind(P^{\otimes m_1})$
or $\ind(P^{\otimes m_1})\nmid \tfrac1{m_1}\ind(P)$. Thus it is enough to prove the claim when $m=p$ is a prime. Write $\ind(P)=p^ac$, where $(p,c)=1$. 
By (\ref{SB.index.thm}) we may assume that $\dim P=p^ac-1$.

By (\ref{|F|.defn})  $P^{\otimes r}$ is similar to $|\o_{P^\vee}(p)|$,
which has dimension $\tbinom{p^ac+p-1}{p}-1$.
Since    $p^{a-1}$ is the largest $p$-power dividing $\tbinom{p^ac+p-1}{p}$,
it is also the largest $p$-power dividing
$\ind(P^{\otimes p})$ by (\ref{SB.index.thm}). \qed

\begin{cor}\label{period.prop.cor} Let $P$ be  a Severi--Brauer variety. Then
$$
\operatorname{period}(P)\mid \ind(P)\mid 
\operatorname{period}(P)^{\dim P}.
$$
In particular, the period and the index have the same prime factors.
\end{cor}

Proof. By (\ref{period.prop}) the period of $P$ divides $\dim P+1$. 
Since $P$ and $P^{\rm min}$  have the same period, 
the period of $P$ divides $\dim P^{\rm min}+1$ and the latter equals
$\ind(P)$ by (\ref{SB.index.thm}).

If $\o_P(m)$ is a line bundle then intersecting the zero set of $\dim P$
general sections yields a 0-cycle of degree $m^{\dim P}$. 
Thus $\ind(P)\mid \operatorname{period}(P)^{\dim P}$.\qed
\medskip

Combining (\ref{period.prop.cor}) with (\ref{submult.index.lem})  gives the following.

\begin{lem}\label{rel.br.dp.lem} Let  $K/k$ be a finite field extension. Then
the kernel of the base-change map  $\br(k)\to \br(K)$ is
$\deg (K/k)$-torsion. \qed
\end{lem}

{\it Aside  \ref{rel.br.dp.lem}.1.} The {\it period-index problem} asks to determine the smallest
number $m$ (depending on the field $k$) such that
$\ind(P)\mid \operatorname{period}(P)^m$
for every Severi--Brauer variety $P$ over $k$. 
See \cite{MR2060023, MR2745688, MR3418522} for various results. 

\begin{prop}[Primary decomposition]  Let $P$ be a minimal 
Severi--Brauer variety. Write $\dim P+1=\prod_i p_i^{c_i}$ where the  $p_i$
are distinct primes. Then there are unique, minimal  Severi--Brauer varieties
$P_i$ such that  $\dim P_i+1= p_i^{c_i}$ and 
$P\cong \otimes_i P_i$.
\end{prop}

Proof. Set $a:=\dim P+1$ and $a_i=p_i^{-c_i}a$. Write  $1=\tsum_i e_ia_i$,
and note that $1\equiv  e_ia_i \mod p_i^{c_i}$.

Assume first that  $P\cong \otimes_i P_i$. By (\ref{period.prop.cor})
 the  period of $P_j$ divides $p_j^{c_j}$, hence
$P_j^{\otimes  e_ia_i}\sim \p^0$ for $j\neq i$. Thus
$$
 P^{\otimes e_{i}a_{i}}\sim \otimes_j P_j^{\otimes  e_ia_i}\sim 
P_i^{\otimes  e_ia_i}\sim P_i,
$$
hence the $P_i$ are unique. Conversely, set  
$P_i=\bigl(P^{\otimes e_{i}a_{i}}\bigr)^{\rm min}$. Then 
$$
P\sim P^{\otimes \sum e_{i}a_{i}} \sim \otimes_i P_i.
$$
Since $P$ is minimal, this shows that
 $P\cong \otimes_i P_i$ iff $\dim P\geq \dim  \otimes_i P_i $.

The period of $P$ divides $a$ by (\ref{period.prop.cor}),
hence the  period of $P_i$ divides $p_i^{c_i}$. Thus the index of $P_i$ is
also a $p_i$-power by (\ref{period.prop.cor}).
On the other hand, by (\ref{ind.of.prods}) the index of $P_i$
 divides $a$, hence the index of $P_i$  divides  $p_i^{c_i}$. Therefore
$\dim P_i\leq p_i^{c_i}-1$ by (\ref{SB.index.thm}.3) hence
$\dim  \otimes_i P_i\leq \dim P$. \qed

\begin{rem} \label{SB.index.thm.rem}
Let $C\subset P$ be a 1-dimensional subscheme. Applying
(\ref{|F|.defn})  to $F=\o_C$ and $m=\pm 1$ to get that 
$
\ind(P)\mid \pm \deg C +\chi(\o_C).
$
This shows that
\begin{enumerate}
\item $ \ind(P)\mid \deg C$ and  $ \ind(P)\mid \chi(\o_C)$ if
$\ind(P)$ is  odd and
\item $ \ind(P)\mid \deg C+\chi(\o_C)$ and  $ \tfrac12\ind(P)\mid \chi(\o_C)$ if
$\ind(P)$ is  even. 
\end{enumerate}
These results are optimal, but  I do not know how to get
necessary and sufficient conditions on the Hilbert polynomials of higher dimensional subvarieties.
\end{rem}

\section{Norm forms and Severi--Brauer varieties}\label{sec.norm.form}

We construct
Severi--Brauer varieties from degree $n$ hypersurfaces in $\p^n$ that are defined using norm forms (\ref{norm.form.defn}) of cyclic Galois extensions.
The main result  is the following.

\begin{thm}\label{normform.SB.thm}
Fix a cyclic Galois extension $K/k$ of degree $n$.
Let $P$ be a  Severi--Brauer variety of dimension $n-1$ over $k$ that has a $K$-point.
 Then $P$ is birational to a  hypersurface of the form
$$
X(K/k, e):=\bigl(\norm_{K/k}(x_1,\dots, x_n)=e t^n\bigr)\subset \p^n
\qtq{for some} e\in k^{\times},
\eqno{(\ref{normform.SB.thm}.1)}
$$
and every such  hypersurface is birational to at least one  Severi--Brauer variety.
\end{thm}

A more precise version is outlined in (\ref{normform.SB.thm.compl}).

\begin{aside} The elements of norm 1 form a $k$-subtorus
$T^{n-1}\subset K^{\times}$. This torus acts on $X(K/k, e)$ such that  
 $X^0(K/k, e):=X(K/k, e)\setminus (t=0)$ is a
torsor under $T^{n-1}$. We will see that  the birational maps
$X(K/k, e)\map P:=P(K/k, e,r)$ are isomorphisms over $X^0(K/k, e)$,
and 
the $P(K/k, e,r)$ inherit this 
$T^{n-1}$-action.  
\end{aside}

\begin{say}[First part of the proof]\label{normform.SB.thm.pf.1}
Let $P$ be a  Severi--Brauer variety of dimension $n-1$ over $k$ that has a $K$-point. Pick a $K$-point $p\in P(K)$ such that its conjugates
$p_1=p, \dots, p_n$ span $P$. By (\ref{crem.tr.exmp.2})   two such points $p, p'$ differ  by an automorphism of $P$, so this is not an additional choice.

Fix also a generator $\sigma \in \gal(K/k)$, this gives a cyclic order on the points $p_i$; we may assume that it is $\sigma(p_i)=p_{i+1}$.
Let 
$$
H_i:=\langle p_1,\dots, p_{i-1}, \widehat{p_i}, p_{i+1},\dots, p_n\rangle
\subset P_K
$$
denote the hyperplane that contains the points $\{p_j:j\neq i\}$ and choose an equation  $H_i=(m_i=0)$ (also over $K$). 
Consider the vectorspace of functions spanned by
$$
 1, \tfrac{m_1}{m_n},\tfrac{m_2}{m_1}, \cdots, \tfrac{m_n}{m_{n-1}}.
$$
Note that although $\sigma(H_i)=H_{i+1}$, it is usually not possible to choose the $m_i$ such that $\sigma^*(m_{i+1})=m_i$ for every $i$. However, no matter how the $m_i$ are chosen, $\sigma^*(m_{i+1})=\lambda_im_i$ for some $\lambda_i\in K^{\times}$, thus the induced linear system on $P$ is defined over $k$.
This gives a rational map
$$
\phi: P\map X\subset \p^n.
\eqno{(\ref{normform.SB.thm.pf.1}.1)}
$$
This may be clearer if we think of $\phi$ as given by the linear system spanned by the hypersurfaces
$$
\begin{array}{l}
H_1+\cdots + H_n \qtq{and}\\
H_1+\cdots + H_{i-2}+ 2H_{i}+H_{i+1}+ \cdots + H_n \qtq{for} i=1,\dots, n.
\end{array}
$$ 
The first of them is defined over $k$, and the others are permuted by $\gal(K/k)$.  So the linear system, and hence $\phi$, are  defined over $k$.

We need to show that $\phi$ is birational and $X$ has the required equation.
Birationality of a map can be checked after a field extension.
Over $K$ we can choose  $u_i:=m_i$ as coordinates, and then the map $\phi_K:\p^{n-1}_{\mathbf u}\cong P_K\map \p^n_{\mathbf x}$ 
is given by
$$
(u_1{:}\dots{:}u_n) \mapsto 
\bigl(1{:}\tfrac{u_1}{u_n}{:}\tfrac{u_2}{u_1}{:} \cdots{:}\tfrac{u_n}{u_{n-1}} 
\bigr),
$$
whose image satisfies the obvious equation
$x_1\cdots x_n=x_0^n$. Over $K$, the inverse is given by
$$
(1{:}x_1{:}\cdots {:} x_n)\mapsto
\bigl(x_1{:}(x_1x_2){:}\cdots {:} (x_1\cdots x_{n})\bigr).
$$
Thus, starting with a  Severi--Brauer variety of dimension $n-1$ over $k$ that has a $K$-point, plus a generator of $\gal(K/k)$, we get a
hypersurface  $X_P$ which is a $k$-form of $(x_1\cdots x_n=x_0^n)$, plus a cyclic 
order on its points of multiplicity $n-1$.

 Next we check that the equation of $X_P$ can be written  in the form (\ref{normform.SB.thm}.1).
\end{say}

\begin{lem}\label{SB.hyper.arith.say.1} Let $k$ be a field and 
$X\subset \p^n$ a hypersurface. Assume that
$$
X_{k^s}\cong (x_1\cdots x_n=t^n).
$$
Then there is a degree $n$ separable $k$-algebra $K$ and $e\in k^{\times}$, such that
$$
X\cong \bigl(\norm_{K/k}({\mathbf x})=et^n\bigr).
$$
Moreover, $K$ is uniquely determined up to isomorphism, and
$e$ is  uniquely determined as an element of
$k^{\times}/\norm_{K/k}(K^{\times})$. 
\end{lem}

Proof.  Note that $(x_1\cdots x_n=t^n)$ has $n$ points of multiplicity $n-1$, and they span a hyperplane. Thus the same holds for $X$.
We can choose $t$ the be the equation of this hyperplane.
Then $X\cap (t=0)$ is a union of   $(n-2)$-planes in $\p^{n-1}\cong (t=0)$, thus its equation is of the form $\bigl(\norm_{K/k}({\mathbf x})=0\bigr)$.  Since $X$ has $n$ points of multiplicity $n-1$, we get the equation
$$
\bigl(\norm_{K/k}({\mathbf x})+L({\mathbf x})t^{n-1}=et^n\bigr),
$$
where $L$ is linear. By assumption $L$ vanishes over $K$, so $L\equiv 0$. \qed

\medskip
Next we write down the construction of $P$ starting from $X_P$.

\begin{say}[Second part of the proof]\label{normform.SB.thm.pf.2}
  Given $\p^{n-1}$ with coordinates  $(u_1{:}\cdots{:}u_n)$, we need to write down its anticanonical embedding using the $x_i=u_i/u_{i-1}$.

  The embedding is given by monomials  $\prod u_i^{c_i}$ such that $c_i\geq 0$ and $\sum c_i=n$. We can scale these by $u_1\cdots u_n$,
  thus work with the  monomials  $\prod u_i^{c_i-1}$.

  The recipe is the following.  Find a string  of the form
  $$
  (c_i, c_{i+1}, \dots, c_{j-1}, c_j)=(0, 1,\dots, 1,>1),
  $$
  divide  $\prod u_i^{c_i-1}$ by  $u_j/u_i=x_jx_{j-1}\cdots x_{i+1}$,
  and repeat. At the end we write $\prod u_i^{c_i-1}$ as a monomial in the $x_j$. The recipe is invariant under cyclic permutations of the $x_i$.
  If $\norm_{K/k}=\prod_i \ell_i$, then replacing the $x_i$ with the
  $\ell_i$ gives a linear system that is Galois invariant. Thus it defines the required inverse of $\phi$ (\ref{normform.SB.thm.pf.1}.1). \qed
  \end{say}

\begin{say}[Complement to (\ref{normform.SB.thm})]\label{normform.SB.thm.compl}
In (\ref{normform.SB.thm})
    the set of such Severi--Brauer varieties forms the group
$\ker[\br(k)\to \br(K)]$, which is $n$-torsion by (\ref{rel.br.dp.lem}).
The isomorphism classes of the  hypersurfaces is determined by
$e\in k^{\times}/\norm_{K/k}(K^{\times})$ by  (\ref{SB.hyper.arith.say.1}).
These two groups are isomorphic, but (\ref{normform.SB.thm}) gives a slightly twisted version of this isomorphism.

The proofs of the following  are left to the reader. Most are straightforward, but (\ref{normform.SB.thm.compl}.7) is cumbersome.
\begin{enumerate} 
\item The construction gives  an isomorphism
$$
\ker[\br(k)\to \br(K)]\times (\z/n\z)^\times\cong 
\bigl(k^{\times}/\norm_{K/k}(K^{\times})\bigr)\times (\z/n\z)^\times.
$$
Thus
for every  hypersurface  as in (\ref{normform.SB.thm}.1) we get birationally equivalent
 Severi--Brauer varieties $P(K/k, e,r)$ for every $r\in (\z/n\z)^{\times}$.
\item  $P(K/k, e,r)$ is trivial over $k$ iff 
$e\in\norm_{K/k}(K^{\times})$, and 
\item $P(K/k, e, r)\sim P(K/k, e, 1)^{\otimes r}$.
\end{enumerate}
The construction of the Severi--Brauer varieties from the
hypersurface using blow-ups and contractions is given by the following steps.
\begin{enumerate}\setcounter{enumi}{3}
\item The divisor $D:= (t=0)\subset X(K/k, e)$ has $n$  geometric irreducible components. Choosing a generator of $\gal(K/k)$  fixes a cyclic order  
$D_1, \dots, D_n$.
\item $X(K/k, e)$ has $n$ geometric points of multiplicity $n-1$.
Let $\pi:Y(K/k, e)\to X(K/k, e)$ denote their blow-up
with exceptional divisor $E$.
\item $E$ has $n(n-1)$ geometric irreducible components $E(i,j)$ naturally indexed by
ordered pairs $\{(i,j): 1\leq i\neq j\leq n\}$, and $n-1$  irreducible components 
$E_r$ such that $E_r=\cup_iE(i,i-r)$.
\item For each  $r\in (\z/n\z)^{\times}$, the image of $Y=Y(K/k, e)$ by the linear system
$$
\bigl|-K_Y+\tfrac{(n+1)(n-2)}{2}\bigl(\pi^*H - E_r\bigr)\bigr|
$$
is an anticanonically embedded  Severi--Brauer variety $P(K/k, e, r)$.
\end{enumerate}

\end{say}

\section{Severi--Brauer schemes}\label{sec.11}

\begin{defn}\label{SBS.defn} A {\it Severi--Brauer scheme} is a 
smooth, proper morphism  $p:P\to S$ all of whose fibers are
Severi--Brauer varieties. Note that $\omega_{P/S}^{-1}$ is $p$-ample, hence 
$p:P\to S$ is projective.
\end{defn}

We aim to show that the basic set-up of Assertion~\ref{ass.1.1}
generalizes to Severi--Brauer schemes. The vector bundle $F(P)$, or rather its dual, appears in  \cite[\S 8.4]{MR0338129}.

\begin{thm}\label{SBS.A.thm} Let  $p:P\to S$ be a  Severi--Brauer scheme.
\begin{enumerate}
\item There is a unique non-split extension
$$
0\to \o_P\to F(P)\to T_{P/S}\to 0
$$
that induces the extension given in (\ref{main.exmp}.1) on all fibers.
\item  $p_*\sEnd\bigl(F(P)\bigr)$ is a locally free sheaf of algebras
whose fiber over $s\in S$ is the central simple algebra corresponding to $P_s$.
\end{enumerate}
\end{thm}

Note that  $p_*\sEnd\bigl(F(P)\bigr)^{\rm opp}$ is called the 
{\it Azumaya algebra} corresponding to $P$;  taking the opposite is just a 
 convention. 
\medskip

Proof. 
We aim to extend the method of Example~\ref{main.exmp} to
Severi--Brauer schemes. 
A priori, the extensions in (\ref{main.exmp}.1) tell us only about the fibers, so  we should look for a non-split extension
$$
0\to p^*L\to F(P)\to T_{P/S}\to 0,
\eqno{(\ref{SBS.A.thm}.3)}
$$
where $L$ is some line bundle on $S$.
The  extension class  is 
$$
\eta\in \ext^1(T_{P/S}, p^*L)=H^1(P, p^*L\otimes \Omega^1_{P/S}).
\eqno{(\ref{SBS.A.thm}.4)}
$$
Note that $p_*\Omega^1_{P/S}=0$, hence
$$
H^1(P, p^*L\otimes \Omega^1_{P/S})=
H^0\bigl(S, L\otimes R^1p_*\Omega^1_{P/S}\bigr).
\eqno{(\ref{SBS.A.thm}.5)}
$$
Since $H^1(P_s, \Omega^1_{P_s})=k$ and $H^i(P_s, \Omega^1_{P_s})=0$ for $i>0$ for every fiber,
$R^1p_*\Omega^1_{P/S}$ is a line bundle. Thus we need to take
$$
L:= \bigl(R^1p_*\Omega^1_{P/S}\bigr)^*.
\eqno{(\ref{SBS.A.thm}.6)}
$$
We check in (\ref{R1.lem}) that
$ R^1p_*\Omega^1_{P/S}\cong \o_S$, hence in fact we do have an extension as in 
(\ref{SBS.A.thm}.1) and the formation of $F(P)$ commutes with arbitrary base change. 

Finally note that $\sEnd\bigl(F(P)\bigr)$ is a locally free sheaf on
$P$ whose restriction to every geometric fiber is a sum of copies of $\o_{P_s}$.
Thus the  formation of  $p_*\sEnd\bigl(F(P)\bigr)$ commutes with arbitrary base change.
\qed

\medskip

In order to understand $R^1\pi_*\Omega^1_{P/S}$, we need a 
better way of writing the Euler sequence.

\begin{say}[Euler sequence, functorially] \label{TPN.2.say}
    For later purposes we work with $\p^n$-bundles.

Let $W$ be a locally free sheaf of rank $n+1$ over $S$, and consider $\pi:P:=\p_S(W)\to S$.
We use  projectivization as in \cite[p.162]{hartsh},  
thus  there are natural identifications
$$
W= \pi_*\o_{P}(1) \qtq{and}
P\cong\bigl(W^*\setminus \{0\}\bigr)/\gm.
\eqno{(\ref{TPN.2.say}.1)}
$$
The cotangent version of the Euler sequence 
is the sheafification of the graded module homomorphism 
$$
e^*:  W\otimes \sym_S(W^*)[-1] \to \sym_S(W^*),  
\qtq{where} e^*(w\otimes w^*):=w(w^*)\in \o_S,
\eqno{(\ref{TPN.2.say}.2)}
$$
 $\sym_S(W)$ denotes the symmetric algebra of $W$ over $S$, and $[-1]$ denotes shifting degrees  by 1. 
Thus the  functorial  cotangent  Euler sequence  is
$$
0\to \Omega_{P/S}^1\to \pi^*W\otimes \o_{P/S}(-1) \stackrel{e^*}{\to} \o_P\to 0.
\eqno{(\ref{TPN.2.say}.3)}
$$
The dual tangent version is
$$
0\to \o_{P}\stackrel{e}{\to}  \pi^*W^*\otimes \o_{P/S}(1)
\stackrel{\partial}{\to}  T_{P/S}\to 0,
\eqno{(\ref{TPN.2.say}.4)}
$$
Both are  unchanged if we tensor $W$ with a line bundle.

 Pushing the latter forward gives 
$$
0\to \o_S\stackrel{e}{\to}  W^*\otimes W
\stackrel{\partial}{\to}  \pi_*T_{P/S}\to 0.
\eqno{(\ref{TPN.2.say}.5)}
$$
\end{say}

\begin{cor} \label{R1.lem} Let $\pi:P\to S$ be a Severi--Brauer scheme.
Then $R^1\pi_*\Omega^1_{P/S}$ is canonically isomorphic to $\o_S$.
\end{cor}

Proof. Assume first that  $\pi:P\to S$ is the projectivization of a vector bundle $W$  of rank $n+1$ over $S$.
Pushing forward  (\ref{TPN.2.say}.3) gives a canonical isomorphism
$$
R^1\pi_*\Omega^1_{P/S}\cong 
\pi_*\o_{P}\cong \o_S,
$$
which is  unchanged if we tensor $W$ with a line bundle.
We are done since a sheaf that has a canonical trivialization \'etale locally is trivial. \qed


\begin{thebibliography}{AHVAV17}

\bibitem[AHVAV17]{MR3587845}
Asher Auel, Brendan Hassett, Anthony V\'{a}rilly-Alvarado, and Bianca Viray
  (eds.), \emph{Brauer groups and obstruction problems}, Progress in
  Mathematics, vol. 320, Birkh\"{a}user/Springer, Cham, 2017, Moduli spaces and
  arithmetic. \MR{3587845}

\bibitem[Art57]{MR0082463}
Emil Artin, \emph{Geometric algebra}, Interscience Publishers, Inc., New
  York-London, 1957. \MR{0082463}

\bibitem[Bae52]{MR0052795}
Reinhold Baer, \emph{Linear algebra and projective geometry}, Academic Press,
  Inc., New York, N.Y., 1952. \MR{0052795}

\bibitem[BLR90]{blr}
Siegfried Bosch, Werner L{\"u}tkebohmert, and Michel Raynaud, \emph{N\'eron
  models}, Ergebnisse der Mathematik und ihrer Grenzgebiete (3), vol.~21,
  Springer-Verlag, Berlin, 1990. \MR{1045822 (91i:14034)}

\bibitem[BN07]{MR2391341}
Indranil Biswas and D.~S. Nagaraj, \emph{Absolutely split real algebraic vector
  bundles over a real form of projective space}, Bull. Sci. Math. \textbf{131}
  (2007), no.~7, 686--696. \MR{2391341}

\bibitem[BN09]{MR2507589}
\bysame, \emph{Vector bundles over a nondegenerate conic}, J. Aust. Math. Soc.
  \textbf{86} (2009), no.~2, 145--154. \MR{2507589}

\bibitem[Bri13]{MR3194649}
Michel Brion, \emph{Homogeneous projective bundles over abelian varieties},
  Algebra Number Theory \textbf{7} (2013), no.~10, 2475--2510. \MR{3194649}

\bibitem[BSS11]{MR2863422}
M.~Blunk, S.~J. Sierra, and S.~Paul Smith, \emph{A derived equivalence for a
  degree 6 del {P}ezzo surface over an arbitrary field}, J. K-Theory \textbf{8}
  (2011), no.~3, 481--492. \MR{2863422}

\bibitem[BSY24]{blanc2024birationalmapsseveribrauersurfaces}
Jérémy Blanc, Julia Schneider, and Egor Yasinsky, \emph{Birational maps of
  {S}everi-{B}rauer surfaces, with applications to {C}remona groups of higher
  rank}, \url{https://arxiv.org/abs/2211.17123}, 2024.

\bibitem[Cam92]{MR1153019}
Peter~J. Cameron, \emph{Projective and polar spaces}, QMW Maths Notes, vol.~13,
  Queen Mary and Westfield College, School of Mathematical Sciences, London,
  1992. \MR{1153019}

\bibitem[Cas06]{MR2264641}
Rey Casse, \emph{Projective geometry: an introduction}, Oxford University
  Press, Oxford, 2006. \MR{2264641}

\bibitem[CH09]{MR2529476}
Jungkai~A. Chen and Christopher~D. Hacon, \emph{On {U}eno's conjecture {K}},
  Math. Ann. \textbf{345} (2009), no.~2, 287--296. \MR{2529476}

\bibitem[Ch{\^a}44]{MR0014720}
Fran{\cedilla{c}}ois Ch{\^a}telet, \emph{Variations sur un th\`eme de {H}.
  {P}oincar\'e}, Ann. Sci. \'Ecole Norm. Sup. (3) \textbf{61} (1944), 249--300.
  \MR{0014720}

\bibitem[Cox64]{MR0172154}
H.~S.~M. Coxeter, \emph{Projective geometry}, Blaisdell Publishing Co. Ginn and
  Co.\, New York-London-Toronto, 1964. \MR{0172154}

\bibitem[CTS21]{MR4304038}
Jean-Louis Colliot-Th\'el\`ene and Alexei~N. Skorobogatov, \emph{The
  {B}rauer-{G}rothendieck group}, Ergebnisse der Mathematik und ihrer
  Grenzgebiete. 3. Folge. A Series of Modern Surveys in Mathematics, vol.~71, Springer, Cham, [2021] \copyright 2021. \MR{4304038}

\bibitem[dJ04]{MR2060023}
A.~J. de~Jong, \emph{The period-index problem for the {B}rauer group of an
  algebraic surface}, Duke Math. J. \textbf{123} (2004), no.~1, 71--94.
  \MR{2060023}

\bibitem[EH16]{eh-3264}
David Eisenbud and Joe Harris, \emph{3264 and all that---a second course in
  algebraic geometry}, Cambridge University Press, Cambridge, 2016.
  \MR{3617981}

\bibitem[Eis95]{eis-ca}
David Eisenbud, \emph{Commutative algebra}, Graduate Texts in Mathematics, vol.
  150, Springer-Verlag, New York, 1995, With a view toward algebraic geometry.
  \MR{1322960 (97a:13001)}

\bibitem[ELW15]{MR3383601}
H\'{e}l\`ene Esnault, Marc Levine, and Olivier Wittenberg, \emph{Index of
  varieties over {H}enselian fields and {E}uler characteristic of coherent
  sheaves}, J. Algebraic Geom. \textbf{24} (2015), no.~4, 693--718.
  \MR{3383601}

\bibitem[Gro68]{gro-bra-I-III}
Alexander Grothendieck, \emph{Le groupe de {B}rauer. {I--III}.}, Dix expos\'es
  sur la cohomologie des sch\'emas, Adv. Stud. Pure Math., vol.~3,
  North-Holland, Amsterdam, 1968, pp.~46--188. \MR{244271}

\bibitem[GSz06]{MR2266528}
Philippe Gille and Tam{\'a}s Szamuely, \emph{Central simple algebras and
  {G}alois cohomology}, Cambridge Studies in Advanced Mathematics, vol. 101,
  Cambridge University Press, Cambridge, 2006. \MR{2266528}

\bibitem[GW10]{MR2675155}
Ulrich G{\"o}rtz and Torsten Wedhorn, \emph{Algebraic geometry {I}}, Advanced
  Lectures in Mathematics, Vieweg + Teubner, Wiesbaden, 2010, Schemes with
  examples and exercises. \MR{2675155}

\bibitem[Har77]{hartsh}
Robin Hartshorne, \emph{Algebraic geometry}, Springer-Verlag, New York, 1977,
  Graduate Texts in Mathematics, No. 52. \MR{0463157 (57 \#3116)}

\bibitem[Her64]{MR0171801}
I.~N. Herstein, \emph{Topics in algebra}, Blaisdell Publishing Co. Ginn and
  Co.\, New York-Toronto-London, 1964. \MR{0171801}

\bibitem[HP47]{hodge-ped}
W.~V.~D. Hodge and D.~Pedoe, \emph{Methods of {A}lgebraic {G}eometry. {V}ols.
  {I--III.}}, Cambridge, at the University Press, 1947. \MR{0028055 (10,396b)}

\bibitem[Jac80]{MR571884}
Nathan Jacobson, \emph{Basic algebra. {II}}, W. H. Freeman and Co., San
  Francisco, Calif., 1980. \MR{571884}

\bibitem[Kle1874]{zbMATH02716823}
Felix Klein, \emph{Nachtrag zu dem zweiten {A}ufsatze {\"u}ber
  {N}icht-{E}uklidische {G}eometrie}, {Math. Ann.} \textbf{7} (1874), 531--537.

\bibitem[KO81]{MR611862}
Michel Kervaire and Manuel Ojanguren (eds.), \emph{Groupe de {B}rauer}, Lecture
  Notes in Mathematics, vol. 844, Springer, Berlin, 1981, Papers from a Seminar
  held at Les Plans-sur-Bex, March 16--22, 1980. \MR{611862}

\bibitem[Kol13]{k-elw}
J{\'a}nos Koll{\'a}r, \emph{{Esnault-Levine-Wittenberg indices}}, December
  2013, \url{https://arxiv.org/abs/1312.3923}.

\bibitem[Kol18]{k-spsb}
\bysame, \emph{Symmetric powers of {S}everi-{B}rauer varieties}, Ann. Fac. Sci.
  Toulouse Math. (6) \textbf{27} (2018), no.~4, 849--862. \MR{3884611}

\bibitem[KS04]{MR2090670}
Daniel Krashen and David~J. Saltman, \emph{Severi-{B}rauer varieties and
  symmetric powers}, Algebraic transformation groups and algebraic varieties,
  Encyclopaedia Math. Sci., vol. 132, Springer, Berlin, 2004, pp.~59--70.
  \MR{2090670 (2005k:14024)}

\bibitem[KSC04]{ksc}
J{\'a}nos Koll{\'a}r, Karen~E. Smith, and Alessio Corti, \emph{Rational and
  nearly rational varieties}, Cambridge Studies in Advanced Mathematics,
  vol.~92, Cambridge University Press, Cambridge, 2004.

\bibitem[Lan02]{lang-alg}
Serge Lang, \emph{Algebra}, third ed., Graduate Texts in Mathematics, vol. 211,
  Springer-Verlag, New York, 2002. \MR{1878556}

\bibitem[Lie07]{MR2309155}
Max Lieblich, \emph{Moduli of twisted sheaves}, Duke Math. J. \textbf{138}
  (2007), no.~1, 23--118. \MR{2309155}

\bibitem[Lie15]{MR3418522}
\bysame, \emph{The period-index problem for fields of transcendence degree 2},
  Ann. of Math. (2) \textbf{182} (2015), no.~2, 391--427. \MR{3418522}

\bibitem[Mil80]{milne}
James~S. Milne, \emph{\'{E}tale cohomology}, Princeton Mathematical Series,
  vol.~33, Princeton University Press, Princeton, N.J., 1980. \MR{559531
  (81j:14002)}

\bibitem[Muk81]{MR607081}
Shigeru Mukai, \emph{Duality between {$D(X)$} and {$D(\hat X)$} with its
  application to {P}icard sheaves}, Nagoya Math. J. \textbf{81} (1981),
  153--175. \MR{607081}

\bibitem[Mum66]{mumf66}
David Mumford, \emph{Lectures on curves on an algebraic surface}, With a
  section by G. M. Bergman. Annals of Mathematics Studies, No. 59, Princeton
  University Press, Princeton, N.J., 1966. \MR{0209285 (35 \#187)}

\bibitem[Nov17]{2017arXiv170103020N}
Sa{\v s}a Novakovi\'c, \emph{{No phantoms in the derived category of curves
  over arbitrary fields, and derived characterizations of Brauer-Severi
  varieties}}, ArXiv e-prints (2017).

\bibitem[Nov24]{MR4796081}
\bysame, \emph{Absolutely split locally free sheaves on proper {$k$}-schemes
  and {B}rauer-{S}everi varieties}, Bull. Sci. Math. \textbf{197} (2024), Paper
  No. 103494, 25. \MR{4796081}

\bibitem[Qui73]{MR0338129}
Daniel Quillen, \emph{Higher algebraic {$K$}-theory. {I}}, Algebraic
  {$K$}-theory, {I}: {H}igher {$K$}-theories ({P}roc. {C}onf., {B}attelle
  {M}emorial {I}nst., {S}eattle, {W}ash., 1972), Springer, Berlin, 1973,
  pp.~85--147. Lecture Notes in Math., Vol. 341. \MR{0338129}

\bibitem[Rey1866]{zbMATH02723115}
Theodor Reye, \emph{Die {G}eometrie der {L}age, 2 vols.}, Hannover,
  R{\"u}mpler, 1866.

\bibitem[Rus1903]{zbMATH02655800}
Bertrand Russell, \emph{The principles of {M}athematics}, Cambridge University
  Press, 1903.

\bibitem[SdJ10]{MR2745688}
Jason Starr and Johan de~Jong, \emph{Almost proper {GIT}-stacks and
  discriminant avoidance}, Doc. Math. \textbf{15} (2010), 957--972.
  \MR{2745688}

\bibitem[Seg42]{MR0008171}
B.~Segre, \emph{The {N}on-singular {C}ubic {S}urfaces}, Oxford University
  Press, Oxford, 1942. \MR{8171}

\bibitem[Ser59]{MR0103191}
Jean-Pierre Serre, \emph{Groupes alg\'ebriques et corps de classes},
  Publications de l'institut de math\'ematique de l'universit\'e de Nancago,
  VII. Hermann, Paris, 1959. \MR{0103191}

\bibitem[Ser79]{MR554237}
\bysame, \emph{Local fields}, Graduate Texts in Mathematics, vol.~67,
  Springer-Verlag, New York-Berlin, 1979, Translated from the French by Marvin
  Jay Greenberg. \MR{554237}

\bibitem[Ser88]{MR918564}
\bysame, \emph{Algebraic groups and class fields}, Graduate Texts in
  Mathematics, vol. 117, Springer-Verlag, New York, 1988, Translated from the
  French. \MR{918564}

\bibitem[Sev1932]{severi-SB}
Francesco Severi, \emph{Un nuovo campo di ricerche nella geometria sopra una
  superficie e sopra una variet{\`a} algebrica}, Mem. Accad. Ital. Mat.
  \textbf{3} (1932), 1--52.

\bibitem[Sko55]{MR73643}
Th. Skolem, \emph{Einige {B}emerkungen \"uber die {A}uffindung der rationalen
  {P}unkte auf gewissen algebraischen {G}ebilden}, Math. Z. \textbf{63} (1955),
  295--312. \MR{73643}

\bibitem[{Sta}22]{stacks-project}
{Stacks Project Authors}, \emph{{S}tacks {P}roject},
  \url{http://stacks.math.columbia.edu}, 2022.

\bibitem[Sza06]{sz-s-b}
Endre Szab{\'o}, \emph{Severi-{B}rauer varieties, (unpublished)}, 2006.

\bibitem[Tre91]{tregub}
S.~L. Tregub, \emph{Birational equivalence of {B}rauer-{S}everi manifolds},
  Uspekhi Mat. Nauk \textbf{46} (1991), no.~6(282), 217--218. \MR{1164209}

\bibitem[vdW91]{vdW-book}
B.~L. van~der Waerden, \emph{Algebra. {V}ol. {I--II}}, Springer-Verlag, New
  York, 1991, Based in part on lectures by E. Artin and E. Noether, Translated
  from the seventh German edition by Fred Blum and John R. Schulenberger.
  \MR{1080172}

\bibitem[Ver82]{MR657419}
Alain H. M.~J. Verschoren (ed.), \emph{Brauer groups in ring theory and
  algebraic geometry}, Lecture Notes in Mathematics, vol. 917, Springer-Verlag,
  Berlin-New York, 1982. \MR{657419}

\bibitem[vS1857]{vStaudt-2}
Karl Georg~Christian von Staudt, \emph{Beitr{\"a}ge zur {G}eometrie der {L}age,
  {H}eft 2}, Kornschen {B}uchhandlung, {N}{\"u}rnberg, 1857.

\bibitem[VY1908]{MR1506049}
Oswald Veblen and John~Wesley Young, \emph{A {S}et of {A}ssumptions for
  {P}rojective {G}eometry}, Amer. J. Math. \textbf{30} (1908), no.~4, 347--380.
  \MR{1506049}

\bibitem[Wei62]{MR0144898}
Andr{\'e} Weil, \emph{Foundations of algebraic geometry}, American Mathematical
  Society, Providence, R.I., 1962. \MR{0144898}

\bibitem[Whi1906]{whitehead}
Alfred~North Whitehead, \emph{The axioms of projective geometry}, Cambridge
  University Press, 1906.

\bibitem[Yos06]{MR2306170}
K{\=o}ta Yoshioka, \emph{Moduli spaces of twisted sheaves on a projective
  variety}, Moduli spaces and arithmetic geometry, Adv. Stud. Pure Math.,
  vol.~45, Math. Soc. Japan, Tokyo, 2006, pp.~1--30. \MR{2306170}

\bibitem[ZS58]{MR0090581}
Oscar Zariski and Pierre Samuel, \emph{Commutative algebra, {V}olume {I}}, The
  University Series in Higher Mathematics, D. Van Nostrand Company, Inc.,
  Princeton, New Jersey, 1958, With the cooperation of I. S. Cohen.
  \MR{0090581}

\end{thebibliography}

\def\cprime{$'$} \def\cprime{$'$} \def\cprime{$'$} \def\cprime{$'$}
  \def\cprime{$'$} \def\dbar{\leavevmode\hbox to 0pt{\hskip.2ex
  \accent"16\hss}d} \def\cprime{$'$} \def\cprime{$'$}
  \def\polhk#1{\setbox0=\hbox{#1}{\ooalign{\hidewidth
  \lower1.5ex\hbox{`}\hidewidth\crcr\unhbox0}}} \def\cprime{$'$}
  \def\cprime{$'$} \def\cprime{$'$} \def\cprime{$'$}
  \def\polhk#1{\setbox0=\hbox{#1}{\ooalign{\hidewidth
  \lower1.5ex\hbox{`}\hidewidth\crcr\unhbox0}}} \def\cdprime{$''$}
  \def\cprime{$'$} \def\cprime{$'$} \def\cprime{$'$} \def\cprime{$'$}
\providecommand{\bysame}{\leavevmode\hbox to3em{\hrulefill}\thinspace}
\providecommand{\MR}{\relax\ifhmode\unskip\space\fi MR }
\providecommand{\MRhref}[2]{%
  \href{http://www.ams.org/mathscinet-getitem?mr=#1}{#2}
}
\providecommand{\href}[2]{#2}

\bigskip

\noindent  Princeton University, Princeton NJ 08544-1000

{\begin{verbatim} kollar@math.princeton.edu\end{verbatim}}

\end{document}